\documentclass[review,10pt]{elsarticle}

\usepackage{etoolbox}
\makeatletter
\patchcmd{\ps@pprintTitle}{\footnotesize\itshape
	Preprint submitted to \ifx\@journal\@empty Elsevier
	\else\@journal\fi\hfill\today}{\relax}{}{}
\makeatother

\usepackage{color,stmaryrd}
\usepackage{natbib}

\usepackage{latexsym}
\usepackage{graphics}
\usepackage{placeins}

\usepackage{amsbsy,textcomp}
\usepackage{amsmath,amsfonts,amssymb}

\definecolor{cyan}{rgb}{0,0.5,0.5}

\usepackage[dvipsnames]{xcolor}
\PassOptionsToPackage{normalem}{ulem}
\usepackage{ulem}
\providecolor{FRC}{rgb}{0,0,1}

\providecolor{SRC}{rgb}{0,0,0}
\newcommand{\SRC}[1]{{\color{SRC}{}#1}}

\journal{}

\usepackage{cite}
\usepackage{amssymb,amsfonts,amsmath}
\usepackage{mathtools} 
\usepackage{color}
\usepackage{graphicx,float}

\usepackage{booktabs}
\usepackage{dcolumn,cancel}
\newcolumntype{d}[1]{D{.}{.}{#1}}  
\usepackage[normalem]{ulem}

\newtheorem{remark}{Remark}

\newcommand{\es}{\hspace*{0.4mm}}

\newcommand{\eps}{\epsilon}

\newcommand{\be}{\boldsymbol{e}}
\newcommand{\bn}{\boldsymbol{n}}

\newcommand{\sip}{\!\cdot\!}
\newcommand{\bI}{\boldsymbol{I}}

\newcommand{\bx}{\boldsymbol{x}}

\newcommand{\by}{\boldsymbol{y}}

\newcommand{\bzero}{\boldsymbol{0}}

\newcommand{\hh}{\hspace*{0.7pt}}

\begin{document}

\begin{frontmatter}

\title{Homogenization of the wave equation with non-uniformly oscillating coefficients}

\author{Danial P. Shahraki}
\author{Bojan B. Guzina \corref{cor1}}
\address{Department of Civil, Environmental and Geo- Engineering, University of Minnesota, Twin Cities, MN 55455, USA}
\cortext[cor1]{Corresponding Author: Bojan B. Guzina; 500 Pillsbury Drive SE, Minneapolis, MN 55455; Email: guzin001@umn.edu; Phone: 612-626-0789.}

\begin{abstract}
The focus of our work is dispersive, second-order effective model describing the low-frequency wave motion in heterogeneous (e.g.~functionally-graded) media endowed with periodic microstructure. For this class of quasi-periodic medium variations, we pursue homogenization of the scalar wave equation in~$\mathbb{R}^d$, $d\geqslant 1$ within the framework of multiple scales expansion. When either $d=1$ or~$d=2$, this model problem bears direct relevance to the description of (anti-plane) shear waves in elastic solids. By adopting the lengthscale of microscopic medium fluctuations as the perturbation parameter, we synthesize the germane low-frequency behavior via a fourth-order differential equation (with smoothly varying coefficients) governing the mean wave motion in the medium, where the effect of microscopic heterogeneities is upscaled by way of the so-called cell functions. In an effort to demonstrate the relevance of our analysis toward solving boundary value problems (deemed to be the ultimate goal of most homogenization studies), we also develop effective boundary conditions, up to the second order of asymptotic approximation, applicable to one-dimensional (1D) shear wave motion in a macroscopically heterogeneous solid with periodic microstructure. We illustrate the analysis numerically in 1D by considering (i) low-frequency wave dispersion, (ii) mean-field homogenized description of the shear waves propagating in a finite domain, and (iii) full-field homogenized description thereof. In contrast to (i) where the overall wave dispersion appears to be fairly well described by the leading-order model, the results in~(ii) and~(iii) demonstrate the critical role that higher-order corrections may have in approximating the actual waveforms in quasi-periodic media. 
\end{abstract}


\end{frontmatter}

\section{Introduction} \label{intro}

\noindent Making use of phenomena such as dispersion, frequency-dependent anisotropy, band gaps, frequency-selective reflection, and negative index of refraction~\citep{Brillouin,Bertoni,Shelby}, phononic materials and periodic composites can be tailored in a way to manipulate waves toward achieving cloaking, vibration control, and sub-wavelength imaging~\citep{Milton, Zhu, Maldovan}. Simulating the underpinning wave motion problem that features rapid (sub-wavelength) variations in the medium, however, can be computationally taxing. To alleviate the impediment, one idea is to ``homogenize'' the medium, i.e. to obtain an effective field equation for the problem that simultaneously: (i)~captures salient attributes of the generated wave motion, and (ii) features ``smoothly-varying'' (typically constant) coefficients that are devoid of rapid variations. 

There is a vast body of literature on the homogenization of wave motion in periodic media. One school of thought -- that is rooted in engineering mechanics and targets an effective description of composites -- is the \emph{Willis' method} of effective constitutive relationships~\citep{Milton, Milton2, Nassar,Meng}. Another keen approach to obtaining the ``macroscopic'' description of periodic media is based on the \emph{Floquet-Bloch} theory~\citep{Floquet, Bloch, Brillouin} which considers the problem eigenfunctions in the form of a plane wave modulated by a periodic function~\citep[e.g.][]{Silva, Ruzzene}. This method is often used to obtain the inherent (multi-valued) dispersion relationship, including band gaps, for periodic media. The third school of thought, rooted in mathematics~\citep{Babuska, Bensoussan, Sanchez, Bakhvalov}, is the method of \emph{multiple scales} expansion where the perturbation parameter is defined as the vanishing ratio between the lengthscale of medium fluctuations and some finite wavelength. 

When multiple-scales homogenization is deployed to describe waves in periodic media, the underpinning asymptotic expansion translates the governing differential equation with rapidly-oscillating coefficients into that with constant coefficients, featuring powers of the perturbation parameter. If the homogenization ansatz is pursued to the leading order only, one obtains the classical homogenization theory and ``effective'' medium properties, for instance the effective elastic moduli and effective mass density in the context of elastodynamics. To extend the dynamic range of such ``long-wavelength'' model, higher-order corrections -- representing singular perturbations of the leading-order effective field equation -- bring about the effects of incipient dispersion and frequency-dependent medium anisotropy. Two-scale homogenization of the wave equation in periodic structures is nowadays well understood. For instance, higher-order asymptotic expansions of the low-frequency (and thus long-wavelength) wave motion in periodic media were studied in one~\citep{Fish2001, Chen, Fish2002a, Bakhvalov2006} and multiple~\citep{Boutin, Fish2002b, Fish2004, Bakhvalov2005, Andrianov2008, Wautier} spatial dimensions. 

To better manipulate waves for the purposes of e.g.~vibration isolation or energy harvesting, however, one may also consider introducing global (i.e.~macroscopic) medium variations that may act \emph{in concert} with their microscopic counterpart. One such example is the concept of~\emph{rainbow trapping} \citep[e.g.][]{Tsak2007,Zhu2013} where an arrangement of dissimilar unit cells -- say sub-wavelength ``resonators'' with progressively varying dynamic characteristics -- is used to enlarge a band gap, relative to what is achievable by purely periodic assemblies. In this setting, it is useful to think of a functionally-graded medium~\citep{Miya2013} endowed with periodic microstructure, see for example the ``staircase'' profiles~\citep{Zhu2013, Colo2016} designed to rainbow trap acoustic and seismic waves. A recent study~\citep{Cell2019} further demonstrates the \emph{spatial arrangement} of a given spectrum of unit cells may significantly affect the performance of a rainbow trap. In the context of optimal design, this exposes the need for rationally constructing a homogenized description of media whose material properties vary on both macroscopic and microscopic scales. As will be seen shortly, our study investigates the low-frequency behavior of such media and does not cover the phenomenon of rainbow trapping; however a generalization of this work in the context of high-frequency homogenization~\citep{Cras2010} could find immediate use in the optimal design of such ``band-stop'' filters.

In the literature, a limited number of works have addressed the homogenization of media with both  macro- and micro-scale variations, especially when considering higher-order asymptotic corrections. Among the earliest studies, one can refer to the multiple-scales frameworks~\citep{Bensoussan, Bakhvalov} that (among other topics) consider elliptic boundary value problems for equations with non-uniformly oscillating coefficients. In~\citep{Andrianov2006}, a novel asymptotic approach was introduced to integrate differential equations governing quasi-periodic structures. The common thread in these developments is their focus on the leading-order homogenized model. Transcending such limitation, a first-order multiscale expansion of the linear elasticity problem for quasi-periodic structures was investigated in~\citep{Cao1999}. More recently, \citep{Fang, Su2011, Ma2013, Yang2018} considered a second-order multiple scales description of the heat conduction, linear elasticity, and thermo-elastic problems in quasi-periodic porous materials. A second-order macroscopic elastic energy of the quasi-periodic materials was investigated in~\citep{le2018second}.

Concerning the wave motion in quasi-periodic media, homogenization of the elastic i.e.~seismic wave equation was pursued via multiple scales expansion up to the leading i.e.~zeroth order for three-dimensional waves~\citep{Cupillard2018}, and up to the first order for two-dimensional (anti-plane shear) SH and (in-plane compressional and shear) P-SV waves~\citep{Guillot,Capdeville}. Recently,~\citep{Dong2017} proposed a second-order, two-scale asymptotic model of the damped wave equation in quasi-periodic media. 

At this point, however, one should note that the first- and second-order models in~\citep{Guillot,Capdeville,Dong2017} are \emph{incomplete} for they disregard asymptotic corrections of the \emph{mean} wave motion. In this vein, the authors in~\citep{Guillot,Capdeville} refer to their expansions as those of ``partial order one''. To highlight the issue, let us denote the featured wave motion by~$u$, its mean-field variation by~$\langle u \rangle$, and the germane perturbation parameter by~$\eps$. In this setting, the multiple-scales approach results in an asymptotic expansion
\[
u(x) \;=\; u_0(x) \,+\, \eps u_1(x,y) \,+\, \eps^2 u_2(x,y) \,+\, \ldots
\] 
where~$x$ and~$y=\eps^{-1}x$ on the right-hand side signify the so-called ``slow'' and ``fast'' variable, respectively, while the asymptotic corrections~$u_j (x,y)$ are built recursively in terms of: (i) lower-order \emph{mean-field} variations~$\langle u_k\rangle(x)$,  $k=\overline{0,j\!-\!1}$ that are governed by the respective (homogenized) field equations, and (ii) fast-oscillating \emph{cell functions} that depend exclusively on medium properties. In periodic media, it is well known~\citep[e.g.][]{Mosk1997,Wautier} that~$\langle u_1\rangle\equiv 0$. In quasi-periodic media, on the other hand, $\langle u_1\rangle\neq 0$ in general. As a result, discarding~$\langle u_j\rangle$, $j\geqslant 1$ as in~\citep{Guillot,Capdeville,Dong2017} leads to incomplete i.e.~``partial'' higher-order solutions.

In this vein, we pursue a (complete) second-order homogenization of the scalar wave equation in~$\mathbb{R}^d$, $d\geqslant 1$ for the class of quasi-periodic media that feature: (i) \emph{smooth} macroscopic variation, and (ii) \emph{periodic} microscopic fluctuation. The analysis commences with a one-dimensional primer ($d\!=\!1$) and demonstrates, via the multiple scales approach, that second-order homogenization of low-frequency wave motion in a medium with non-uniformly (yet ``rapidly'') oscillating coefficients yields a fourth-order effective differential equation with smoothly varying coefficients. This result is then extended to describe the effective wave motion in quasi-periodic media for $d>1$. Motivated by a recent (one-dimensional) study of the effective wave motion in \emph{bounded periodic domains}~\citep{Corna2019}, we next develop effective ``Dirichlet'' and ``Neumann'' boundary conditions (up to the second order of asymptotic approximation) for one-dimensional waves in quasi-periodic media, and we use this result to homogenize one-dimensional \emph{boundary value problem} for a quasi-periodic domain of finite extent. A set of numerical results including (i) dispersion curves, (ii) mean-filed approximations, and (iii) full-field approximations of the wave motion in quasi-periodic media is included to illustrate the utility of the proposed homogenization framework. In contrast to (i) where the overall wave dispersion (due to combined ``action'' of  micro- and macro-scale heterogeneities) appears to be well described by the leading-order model, the results in~(ii) and~(iii) highlight the critical role that  higher-order corrections have in maintaining the fidelity of a homogenized description. To the authors' knowledge, this is the first study where the effective boundary conditions in quasi-periodic media have been considered.

\section{Problem statement} \label{ScWv}

\begin{figure}[h!]  \vspace*{-1mm}
\centering{\includegraphics[width=0.55\linewidth]{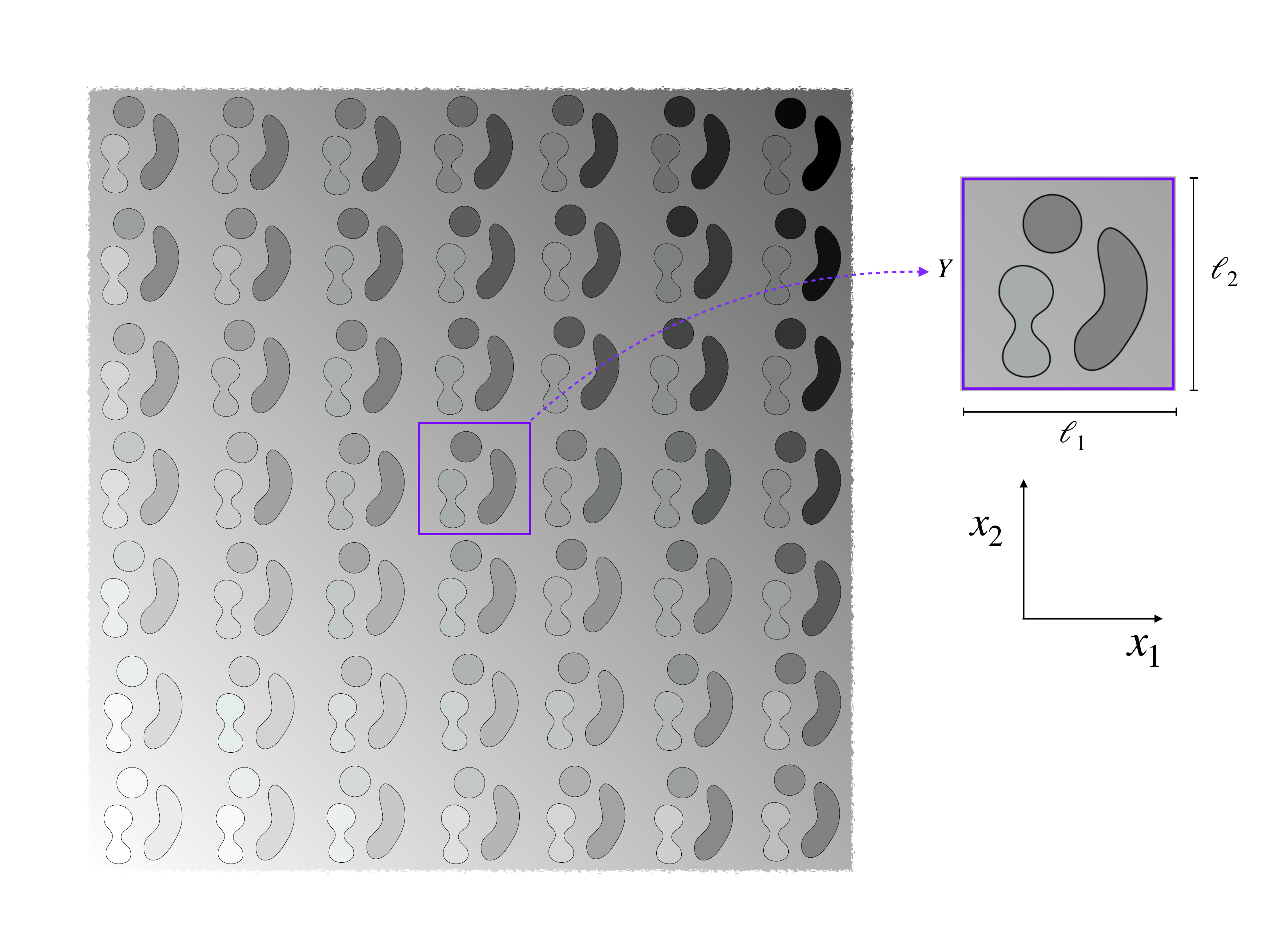}}\vspace*{-2mm}
\caption{Quasi-periodic medium featuring (i) smooth macroscopic variation, and (ii) periodic (but not necessarily smooth) microscopic fluctuation.} 
\label{Sche_HetMat_2d}
\end{figure}

\noindent Assuming all parameters and variables hereon to be dimensionless (with reference to a suitable dimensional platform), we consider the scalar wave motion in an unbounded heterogeneous medium, namely 
\begin{eqnarray} \label{WaveEq}
\nabla \!\cdot\! \big(G (\bx ) \nabla u \big) + \rho(\bx)\es \omega^2 u \;=\; 0, \qquad \bx\in\mathbb{R}^d, 
\end{eqnarray}
where~$d\geqslant 1$ and ~$\omega$ denotes the oscillation frequency. With reference to Fig.~\ref{Sche_HetMat_2d}, we let the propagation medium feature both \emph{smooth} (but otherwise arbitrary) \emph{macroscopic variation}, and \emph{periodic small-scale fluctuation}; specifically, we assume the coefficients~$G(\bx)>0$ and~$\rho(\bx)>0$ in~\eqref{WaveEq} to admit either additive or multiplicative separation between ``macroscopic'' and ``microscopic'' variations according to  
\begin{alignat}{3}
\text{additive separation:} && ~~ G(\bx) = G'(\bx)+G''(\bx/\epsilon), \quad 
&& \rho(\bx) = \rho'(\bx)+\rho''(\bx/\epsilon) \label{HetMat1}\\
\text{multiplicative separation:} && G(\bx)= G'(\bx)\, G''(\bx/\epsilon), \quad && \rho(\bx) = \rho'(\bx)\, \rho''(\bx/\epsilon) \label{HetMat2}
\end{alignat}
where $\epsilon=o(1)$ is the germane perturbation parameter; $G'>0$ and $\rho'>0$ are bounded and smooth functions supported in~$\mathbb{R}^d$; and $G''$ and $\rho''$ are $Y$-periodic functions, $Y\subset\mathbb{R}^d$. Making reference to a Cartesian coordinate system tied to an orthonormal basis $\be_j (j=\overline{1,d})$, we define the unit cell of periodicity as  
\begin{eqnarray} \label{MicroS}
Y =\{\bx:0<\bx\!\cdot\!\boldsymbol{e}_j<\ell_j, ~ j =\overline{1,d}\}, \qquad |Y|=1. 
\end{eqnarray}

\begin{remark}
In the context of linear elasticity and anti-plane shear waves, $u, G$ and~$\rho$ in~\eqref{WaveEq} can be interpreted (for~$d\in\{1,2\}$) as the transverse displacement, shear modulus, and mass density, respectively. 
\end{remark}

\begin{remark}
In situations where~$G''(\bx/\eps)$ and~$\rho''(\bx/\eps)$ are discontinuous, (\ref{WaveEq}) is implicitly complemented by the ``perfect bonding'' conditions between smooth (microscopic) constituents. Letting~$\Gamma_{\!\eps}\subset\mathbb{R}^{d-1}$ denote the union of all such material discontinuities in~$\mathbb{R}^d$, we hereon implicitly assume that  
\begin{eqnarray} \label{continuity}
\llbracket  u \rrbracket = 0, ~~~  \llbracket  \boldsymbol{n} \!\cdot\! (G  \nabla u ) \rrbracket =0, \qquad\bx\in\Gamma_{\!\eps}
\end{eqnarray}
where~$\bn$ denotes the unit normal on~$\Gamma_{\!\eps}$, and 
\begin{eqnarray} \label{Jum}
\llbracket  g \rrbracket(\bx) = \lim_{\eta \to 0^+} \big\{g(\bx+\eta \boldsymbol{n}(\bx))- g(\bx-\eta \boldsymbol{n}(\bx))\big\}, \qquad \bx\in\Gamma_{\!\eps} \notag
\end{eqnarray}
signifies the jump across the interface. 
\end{remark}

\subsection{Objective and assumptions}\label{OandA}

\noindent At this point, we observe from the representation~\eqref{HetMat1}--\eqref{HetMat2} of~$G(\bx)$ and~$\rho(\bx)$ that the ``microscopic'' coefficient fluctuations in~\eqref{WaveEq} are in fact $\eps Y$-periodic. In this setting, our goal is to obtain an effective i.e.~macroscopic representation of~\eqref{WaveEq}, given by a field equation (with smoothly varying coefficients) governing the ``\emph{mean}'' wave motion, assuming that:
\begin{itemize}
\item the frequency of oscillation is finite, namely $\omega=O(1)$, yet sufficiently low so that the germane wave dispersion -- as driven by the ``macro'' and ``micro'' medium heterogeneities -- resides inside the \emph{(apparent) first pass band}~\citep{Brillouin}; and  

\item the characteristic lengthscale of~$\eps Y$ is much smaller than the $\omega$-induced dominant wavelength of~$u(\bx)$, as implied by the premise~$\eps = o(1)$. 
\end{itemize}

Motivated by the earlier studies of waves in periodic media~\citep{Boutin, Fish2002b, Fish2004, Bakhvalov2005, Andrianov2008}, we tackle the problem via multiple-scales asymptotic expansion~\citep{Babuska, Bensoussan} that revolves around the mapping $u(\bx)\mapsto u(\bx,\bx/\eps)$ as a means to parse the macroscopic and microscopic wavefield fluctuations. In this setting, we also note that the foregoing restriction on~$\omega$ amounts to considering the  low-frequency, low-wavenumber (LF-LW) homogenization~\citep{Meng,guzina2019rational} of wave motion. For completeness we pursue the asymptotic approximation \emph{up to the second order}, which brings about an effective field equation that captures with high fidelity the combined effects of macroscopic and microscopic medium heterogeneities. 

\begin{remark}
In periodic media, the sought second-order correction -- appearing as a singular perturbation of the leading-order effective field equation~\citep{Andrianov2008, Wautier} -- has only a moderate effect~\citep{Wautier} on the (frequency-dependent) phase of the solution. However, a recent study~\citep{Corna2019} demonstrates that such correction may have a major effect of the solution amplitude when considering wave motion in bounded micro-structured domains. 
\end{remark}

In what follows, we pursue the LF-LW asymptotic treatment i.e.~homogenization of~\eqref{WaveEq}, assuming either~\eqref{HetMat1} or~\eqref{HetMat2}, in terms of the perturbation parameter $\eps$. For brevity of presentation, we describe in detail the homogenization of one-dimensional (1D) waves, and follow up by presenting only the final result for the general case in~$\mathbb{R}^d$, $d\geqslant 1$ -- obtained in an analogous fashion. We then complete the 1D analysis by obtaining (also via multiple scales expansion) the effective boundary conditions applicable to the mean-field motion, and we illustrate numerically the analytical developments by considering 1D waves in both unbounded and bounded quasi-periodic domains.

\section{Effective field equation for one-dimensional problems} \label{MuSc}

\noindent When~$d=1$, field equation~\eqref{WaveEq} reduces to 
\begin{eqnarray} \label{1D_Wave}
\frac{\text{d} }{\text{d} x} \Big( G(x) \frac{\text{d} u}{\text{d} x} \Big) +\rho(x)\es \omega^2 u = 0, \qquad x\in \mathbb{R}. 
\end{eqnarray}
Here~$G(x)$ and~$\rho(x)$ satisfy the one-dimensional counterparts of either~\eqref{HetMat1} or~\eqref{HetMat2}, where~$G''(x/\eps)$ and~$\rho''(x/\eps)$ are $\eps Y$-periodic with~$Y=(0,1)$. Within the framework of two-scale homogenization~\citep{Bensoussan}, this motivates introduction of the ``fast'' coordinate $y=\eps^{-1} x$ and affiliated mappings
\begin{eqnarray} \label{sub}
\begin{split}
& u(x) \mapsto u(x,y),\qquad  \frac{\text{d} }{\text{d} x}  \mapsto \frac{\partial}{\partial x} + \eps^{-1} \frac{\partial}{\partial y},  \\
&  G(x) \mapsto G(x,y), \qquad \rho(x) \mapsto \, \rho(x,y),
\end{split}
\end{eqnarray}
designed to separate the \emph{macroscopic} variations (described in terms of~$x$) from their \emph{microscopic} counterparts (by definition $Y$-periodic in terms of~$y$). In terms of~\eqref{HetMat1}--\eqref{HetMat2}, we specifically see that 
\begin{alignat}{3}
\text{additive separation:} && ~~ G(x,y) = G'(x)+G''(y), \quad 
&& \rho(x,y) = \rho'(x)+\rho''(y) \label{HetMat11d}\\
\text{multiplicative separation:} && G(x,y)= G'(x)\, G''(y), \quad && \rho(x,y) = \rho'(x)\, \rho''(y) \label{HetMat21d}
\end{alignat}
For convenience, we also introduce a two-scale flux quantity~$\sigma$ -- namely the shear stress in terms of anti-plane shear waves -- and its mapping via  
\begin{eqnarray} \label{stress1}
\sigma(x) := G(x) \frac{\text{d}u }{\text{d} x} ~~\mapsto~~ \sigma(x,y) = G(x,y)\Big(\frac{\partial u}{\partial x} + \eps^{-1} \frac{\partial u}{\partial y}\Big) . 
\end{eqnarray}

On substituting~(\ref{sub}) into~(\ref{1D_Wave}) and using short-hand notation $(\cdot)_{,x}=\partial(\cdot)/\partial x$ and $(\cdot)_{,y}=\partial(\cdot)/\partial y$, (\ref{1D_Wave}) can be rewritten in powers of $\epsilon$ as
\begin{eqnarray} \label{ReScEq}
\epsilon^{-2}\big([G u_{, y}]_{,y}\big) + \epsilon^{-1}\big([G  u_{,y}]_{,x} + [G u_{,x}]_{,y}\big) +\big(\rho\es \omega^2 u + [G u_{,x}]_{,x}\big) =0, \qquad x\in\mathbb{R}, ~ y\in Y. 
\end{eqnarray}
We next pursue the asymptotic solution of~(\ref{ReScEq}) in terms of the ansatz  
\begin{eqnarray} \label{AsEx}
u(x,y)=\sum_{k=0}^\infty \epsilon^k u_k (x,y), \qquad \sigma(x,y)=\sum_{k=0}^\infty \epsilon^k \sigma_k (x,y),
\end{eqnarray}
where $\sigma_k = G(u_{k,x} + u_{k+1,y})$. With such definitions, the homogenized i.e. macroscopic LF-LW description of the problem can be effected in terms of the ``mean'' wave motion 
\begin{eqnarray} \label{mean}
\langle u \rangle(x) \,:=\; \int_{0}^{1}u(x,y) \, \text{d}y \;=\; \sum_{k=0}^\infty \eps^k \langle u_k \rangle(x).
\end{eqnarray}
In the sequel, we shall also make use of the partial sums 
\begin{eqnarray} \label{resEx}
u^{[p]}(x,y)=\sum_{k=0}^{p} \epsilon^k u_k (x,y), \qquad \sigma^{[p]}(x,y)=\sum_{k=0}^{p} \epsilon^k \sigma_k (x,y), \qquad \langle u \rangle^{[p]}(x)  \;=\; \sum_{k=0}^{p} \eps^k \langle u_k \rangle(x).
\end{eqnarray}

\subsection{O(1) homogenization} \label{0th_H}

\noindent By virtue of~\eqref{AsEx}, the $O(\epsilon^{-2})$ statement of~\eqref{ReScEq} becomes 
\begin{eqnarray} \label{eps_m2}
[G u_{0,y}]_{,y} =0 \quad \Rightarrow \quad u_{0,y} = G^{-1} c(x), \qquad y\in Y 
\end{eqnarray}
\noindent where $x\in\mathbb{R}$ hereon, and $c(x)$ is a constant of integration. On averaging the last result, we obtain 
\begin{eqnarray} \label{u0}
\langle u_{0,y}\rangle \;=\; c(x) \langle G^{-1}\rangle. 
\end{eqnarray}
Thanks to the $Y$-periodicity of~$u_0(\cdot,y)$, the left-hand side of (\ref{u0}) vanishes identically whereby~$c(x)=0$. Accordingly, the leading-order solution depends exclusively on the macroscopic variable and we write
\begin{eqnarray} \label{u0}
u_0(x,y) \,=\, \langle u_0\rangle. 
\end{eqnarray}
On collecting the $O(\epsilon^{-1})$ terms in~\eqref{ReScEq}, one arrives at 
\begin{eqnarray} \label{eps_m1}
[G (u_{1,y}+ u_{0,x}) ]_{,y} =0, \qquad y\in Y.
\end{eqnarray}

\noindent Since~\eqref{eps_m1} is a linear ordinary differential equation in~$y$, its general solution can be written as 
\begin{eqnarray} \label{u1}
u_1(x,y)  \,=\, \langle u_1\rangle + P(x,y)\es \langle u_0\rangle_{\!,x}, \qquad \langle P \rangle =0, 
\end{eqnarray}

\noindent where $P$ is a zero-mean ($Y$-periodic in~$y$) cell function that satisfies
\begin{eqnarray} \label{P}
[G(1+P_{\!,y})]_{,y}=0, \qquad y\in Y
\end{eqnarray}
for each~$x\in\mathbb{R}$. Next, we proceed to the $O(1)$ governing equation which reads
\begin{eqnarray} \label{eps_0}
[G (u_{2,y} + u_{1,x})]_{,y} + [G (u_{1,y}+u_{0,x}) ]_{,x} + \rho\es \omega^2 \langle u_0\rangle \;=\;0, \qquad y\in Y. 
\end{eqnarray}
Again, (\ref{eps_0}) is a linear equation whereby its solution admits the representation 
\begin{eqnarray} \label{u2}
\begin{split}
u_2(x,y) \,=\, \langle u_2\rangle + P(x,y) &  \langle u_1\rangle_{\!,x} +  \tilde{P}(x,y)  \langle u_0\rangle_{\!,x} + Q(x,y)  \langle u_0\rangle_{\!,xx}, \qquad 
\langle \tilde{P}\rangle = \langle Q\rangle =0,
\end{split}
\end{eqnarray}
where~$\tilde{P}$ and~$Q$ are zero-mean cell functions satisfying 
\begin{eqnarray} \label{Ptilde_Q}
\Big([G (1+P_{\!,y})]_{,x} + [G (P_{,x}+\tilde{P}_{,y})]_{,y}\Big) \langle u_0\rangle_{\!,x}  +\Big(G (1+P_{\!,y}) +[G (P+Q_{,y})]_{,y} \Big) \langle u_0\rangle_{\!,xx} + \rho \es \omega^2 \langle u_0\rangle \;=\; 0, \quad y\in Y
\end{eqnarray}
thanks to~\eqref{u1}. 

\begin{remark}
Up to this point, the key structural difference between the present problem and that describing the wave motion in periodic media~\citep{Wautier} is the presence of the cell function~$\tilde{P}(x,y)$ in the expression~\eqref{u2} for~$u_2$. As will be shown in the sequel, this function vanishes identically when $G'(x)$ and $\rho'(x)$ in~\eqref{HetMat11d}--\eqref{HetMat21d} are constant. 
\end{remark} 

Integrating (\ref{eps_0}) over~$Y$, taking into account the inherent $Y$-periodicity in the fast variable, and assuming $P,Q$ and $\tilde{P}$ to be bounded, we obtain the (leading-order) homogenized field equation governing~$\langle u_0\rangle$, namely 
\begin{eqnarray} \label{O_1_Eq}
[\mu^{\mbox{\tiny{(0)}}}(x) \langle u_0\rangle_{\!,x}]_{,x} + \varrho^{\mbox{\tiny{(0)}}}(x) \omega^2 \langle u_0\rangle  = 0, \qquad x\in\mathbb{R}
\end{eqnarray}
\noindent whose \emph{macroscopically-heterogeneous} effective coefficients are given by 
\begin{eqnarray} \label{rho0_mu0}
\mu^{\mbox{\tiny{(0)}}}(x)= \langle G(1+P_{\!,y}) \rangle, \qquad \varrho^{\mbox{\tiny{(0)}}}(x)=\langle \rho \rangle.
\end{eqnarray}
In contrast to the purely periodic case~\citep{Fish2001,Wautier}, we note that even the leading-order effective model is in this case \emph{dispersive} due to macroscopic variation~\eqref{rho0_mu0} of the effective medium properties. We also note that exposing the (heterogeneous) effective medium properties requires knowledge of the cell function~$P(x,y)$ whose evaluation, and that of its companions, is addressed in Section~\ref{Cell_Func}. 

\subsection{ $O(\epsilon)$ homogenization} \label{first_H}

\noindent On collecting the $O(\epsilon)$ terms in~(\ref{ReScEq}), we find 
\begin{eqnarray} \label{eps_1}
[G (u_{3,y}+u_{2,x}) ]_{,y} +  [G (u_{2,y}+u_{1,x} ) ]_{,x} + \rho\es \omega^2 u_1 \;=\; 0, \qquad y\in Y. 
\end{eqnarray}
By analogy to earlier treatment, the linearity of~\eqref{eps_1} allows us to write the general solution as 
\begin{multline} \label{u3}
u_3(x,y) \,=\, \langle u_3\rangle + P(x,y) \langle u_2\rangle_{\!,x} + \tilde{P}(x,y) \langle u_1\rangle_{\!,x} + Q(x,y) \langle u_1\rangle_{\!,xx} \, +  \\
\tilde{R}(x,y) \langle u_0\rangle_{\!,x} + \tilde{Q}(x,y) \langle u_0\rangle_{\!,xx} + R(x,y) \langle u_0\rangle_{\!,xxx}, \qquad
\langle\tilde{Q}\rangle = \langle R\rangle = \langle\tilde{R}\rangle = 0, 
\end{multline}
\noindent where the zero-mean cell functions~$\tilde{Q}$, $R$ and $\tilde{R}$ are such that 
\begin{multline} \label{WQR}
\big([G(1+P_{\!,y})]_{,x} + [G(\tilde{P}_{,y}+P_{,x})]_{,y}\big) \langle u_1\rangle_{\!,x} + 
\big([G(Q_{,y}+P)]_{,y} + G(1+P_{\!,y}) \big) \langle u_1\rangle_{\!,xx}+  \\
\big([G(P_{,x}+\tilde{P}_{,y})]_{,x} + [G(\tilde{R}_{,y}+\tilde{P}_{,x})]_{,y} +\rho  ~\omega^2P \big) \langle u_0\rangle_{\!,x}+ 
\big([G(\tilde{Q}_{,y}+Q_{,x}+\tilde{P})]_{,y} + [G(P+Q_{,y})]_{,x}+G(P_{,x}+\tilde{P}_{,y})  \big) \langle u_0\rangle_{\!,xx}  \\
 + \big([G(R_{,y}+Q)]_{,y} + G(P+Q_{,y})  \big) \langle u_0\rangle_{\!,xxx}+\rho~\omega^2 \langle u_1\rangle  =0,
\end{multline}

\noindent On integrating (\ref{WQR}) over~$Y$ and exploiting the periodicity of featured quantities, we obtain the effective field equation governing the first-order corrector $\langle u_1\rangle$ as
\begin{eqnarray} \label{O_eps_Eq}
\begin{split}
[\mu^{\mbox{\tiny{(0)}}}(x) \langle u_1\rangle_{\!,x}]_{,x} +\varrho^{\mbox{\tiny{(0)}}}(x) \omega^2 \langle u_1\rangle \;=\;  -[\eta_{,x}(x)+\varrho^{\mbox{\tiny{(1)}}}(x) \omega^2] \langle  u_0\rangle_{\!,x} - [\eta(x)+\mu^{\mbox{\tiny{(1)}}}_{,x}(x)]\langle u_0\rangle_{\!,xx} - \mu^{\mbox{\tiny{(1)}}}(x) \langle u_0\rangle_{\!,xxx},
\end{split}
\end{eqnarray}
\noindent whose ``functionally-graded'' coefficients are given by 
\begin{eqnarray} \label{xiC}
\mu^{\mbox{\tiny{(1)}}}(x) = \langle  G(P+Q_{,y}) \rangle, \qquad \varrho^{\mbox{\tiny{(1)}}}(x) = \langle  \rho P \rangle, \qquad 
\eta(x) = \langle  G(P_{,x}+\tilde{P}_{,y}) \rangle. 
\end{eqnarray}

\subsection{$O(\epsilon^2)$ homogenization} \label{second_H}

\noindent On collecting the $O(\epsilon^2)$ terms in (\ref{ReScEq}), one finds that 
\begin{eqnarray} \label{eps_2}
[G ( u_{4,y}+ u_{3,x}) ]_{,y} +  [G (u_{3,y}+u_{2,x}) \big ]_{,x} + \rho\es\omega^2 u_2 \;=\; 0,
\end{eqnarray}
Following the previously established analysis and integrating~\eqref{eps_2} over $Y$, we obtain the effective equation for~$\langle u_2\rangle$, namely
\begin{multline} \label{O_eps2_Eq}
[\mu^{\mbox{\tiny{(0)}}}(x) \langle u_2\rangle_{\!,x}]_{,x} +\varrho^{\mbox{\tiny{(0)}}}(x)\es \omega^2 \langle u_2\rangle  \;=\; 
-[\eta _{,x}(x) + \varrho^{\mbox{\tiny{(1)}}}(x)\es \omega^2]\langle u_1\rangle_{\!,x} - [\eta (x) + \mu^{\mbox{\tiny{(1)}}}_{,x}(x)] \langle u_1\rangle_{\!,xx} - \mu^{\mbox{\tiny{(1)}}}(x) \langle u_1\rangle_{\!,xxx}  \\
~~-[\phi_{,x}(x) + \tilde{\varrho}^{\mbox{\tiny{(2)}}}(x) \omega^2]\langle u_0\rangle_{\!,x} - [\psi_{,x}(x)+\phi(x)+\varrho^{\mbox{\tiny{(2)}}}(x) \omega^2]\langle u_0\rangle_{\!,xx} 
-[\psi(x)+\mu^{\mbox{\tiny{(2)}}}_{,x}(x)]\langle u_0\rangle_{\!,xxx} - \mu^{\mbox{\tiny{(2)}}} (x) \langle u_0\rangle_{\!,xxxx}
\end{multline}
\noindent where
\begin{eqnarray} \label{phiC}
\begin{split}
& \mu^{\mbox{\tiny{(2)}}}(x) = \langle  G (Q+R_{,y})\rangle, \qquad  \varrho^{\mbox{\tiny{(2)}}}(x) = \langle  \rho Q\rangle, \qquad \tilde{\varrho}^{\mbox{\tiny{(2)}}}(x) = \langle  \rho \tilde{P}\rangle, \\
& \phi(x) = \langle  G (\tilde{P}_{,x}+ \tilde{R}_{,y} ) \rangle, \qquad \psi(x) = \langle G(Q_{,x}+\tilde{P}+ \tilde{Q}_{,y}) \rangle.
\end{split}
\end{eqnarray}

\subsection{Macroscopic description of the mean wave motion} \label{MeanMotion}

\noindent By virtue of~\eqref{mean}, one my conveniently compute the weighted sum $\eqref{O_1_Eq}+\eps\eqref{O_eps_Eq}+\eps^2\eqref{O_eps2_Eq}$, resulting in the second-order effective equation 
\begin{eqnarray} \label{MMEq}
\begin{split}
&\mu^{\mbox{\tiny{(0)}}}(x) \es\langle u\rangle_{\!,xx} +\mu^{\mbox{\tiny{(0)}}}_{,x}(x)\es \langle u\rangle_{\!,x} +\varrho^{\mbox{\tiny{(0)}}}(x) \es\omega^2 \langle u\rangle = \\
&\qquad - \epsilon \Big(  [\eta _{,x}(x) + \varrho^{\mbox{\tiny{(1)}}} (x) \omega^2 ]\langle u\rangle_{\!,x} + [\eta  (x) + \mu^{\mbox{\tiny{(1)}}}_{,x} (x) ] \langle u\rangle_{\!,xx} +  \mu^{\mbox{\tiny{(1)}}} (x) \langle u\rangle_{\!,xxx}  \Big)  \\
&\qquad - \epsilon^2  \Big( [\phi_{,x} (x)+ \tilde{\varrho}^{\mbox{\tiny{(2)}}} (x) \omega^2 ] \langle u\rangle_{\!,x}  +     [  \psi_{,x}(x) +\phi (x)  + \varrho^{\mbox{\tiny{(2)}}} (x)\omega^2 ] \langle u\rangle_{\!,xx} \\
&\qquad\qquad~~  + [\psi (x)+ \mu^{\mbox{\tiny{(2)}}}_{,x}(x)] \langle u\rangle_{\!,xxx} + \mu^{\mbox{\tiny{(2)}}} (x) \langle u\rangle_{\!,xxxx} \Big) + O(\epsilon^3), 
\end{split}
\end{eqnarray}
which demonstrates that an $O(\eps)$ (resp. $O(\eps^2)$) correction of the leading-order effective model 
\[
[\mu^{\mbox{\tiny{(0)}}}(x) \langle u\rangle_{\!,x}]_{,x} +\varrho^{\mbox{\tiny{(0)}}}(x)\es \omega^2 \langle u\rangle \;=\; 0
\]
due to presence of two-scale medium variations entails spatial derivatives up to order three (resp. four).  One key difference between the current model and its counterparts for periodic media~\citep{Fish2001,Wautier} where $G(x,y)=G(y)$ and~$\rho(x,y)=\rho(y)$, however, is that~\eqref{MMEq} features a \emph{nontrivial $O(\eps)$ correction} that can be shown to vanish identically in the periodic case~\citep{guzina2019rational}. With reference to~\eqref{HetMat11d}, for instance, it can be specifically shown that 
\[
\begin{split}
& \bullet~ G'(x)=\rho'(x)=\text{const.} \text{ (no macro-scale heterogeneities)} \\
& \qquad \Longrightarrow \quad  \mu^{\mbox{\tiny{($k$)}}}=\text{const.}, \quad \rho^{\mbox{\tiny{($k$)}}}=\text{const.} ~~(k=0,2), \quad \mu^{\mbox{\tiny{(1)}}}=\varrho^{\mbox{\tiny{(1)}}} = \tilde{\varrho}^{\mbox{\tiny{(2)}}} =\eta = \phi = \psi = 0; \\
& \bullet~ G''(y)=\rho''(y)= 0 \text{ (no micro-scale heterogeneities)} \\
& \qquad \Longrightarrow \quad \mu^{\mbox{\tiny{(0)}}}(x) = G'(x), \quad \varrho^{\mbox{\tiny{(0)}}}(x) = \rho'(x), \quad 
 \mu^{\mbox{\tiny{($k$)}}}=\rho^{\mbox{\tiny{($k$)}}}= \tilde{\varrho}^{\mbox{\tiny{(2)}}}= \eta = \phi = \psi = 0; ~~(k\geqslant 1).
\end{split}
\]
For completeness, we also recall that the effective coefficients $\varrho^{\mbox{\tiny{($k$)}}}$, $\mu^{\mbox{\tiny{($k$)}}} ~ (k = 0,1,\ldots)$, $\eta, \phi$, and $\psi$ are specified via~\eqref{rho0_mu0}, \eqref{xiC} and~\eqref{phiC} in terms of the respective cell functions whose evaluation is examined next.

\begin{remark}
In the context of apparent complexity characterizing~\eqref{MMEq}, a fair question to ask concerns the utility of homogenized descriptions for this class of problems: what is to be gained? Assuming multi-dimensional problems (a subject that will be addressed shortly), a short answer is ``the computational efficiency'' for such homogenized models are characterized by smooth coefficient variations. This is in contrast to the original quasi-periodic medium that may feature material discontinuities at the microscopic ($O(\eps)$) scale, and thus necessitate spatial discretization whose ``fine'' lengthscale is~$o(\eps)$, i.e. several decades smaller than the dominant wavelength.
\end{remark}

\subsection{Cell functions} \label{Cell_Func}

\noindent For generality, we assume (as examined earlier) that the scaled unit cell of periodicity~$Y=(0,1)$ is composed of~$N$ smoothly heterogeneous pieces i.e. sub-cells~$Y_q$, $q=\overline{1,N}$ such that~$G''(y)$ and~$\rho''(y)$ are differentiable within each (open) set~$Y_q$. In this setting, by recalling~\eqref{P} we find that the cell function~$P(x,y)$ featured in~\eqref{u1} solves the boundary value problem 
\begin{eqnarray} \label{CF_P}
\begin{split}
& [G(1+P_{\!,y})]_{,y}=0,\quad y \in  Y_q \\*[-1mm]
& P, ~ G(1+P_{\!,y}) \quad Y\text{-periodic}; \quad  \langle  P \rangle = 0  \\*[-1mm]
& \, \llbracket P\rrbracket =0, \quad \llbracket G(1+P_{\!,y})\rrbracket =0, ~\quad y \in \partial Y_q\backslash \partial Y   
\end{split}
\end{eqnarray}
for any given~$x\in\mathbb{R}$, where the term~$G(1+P_{\!,y})$ can be understood as a generalized flux relevant to~$P$. 

From the leading-order field equation~\eqref{O_1_Eq}, $\langle u_0\rangle$ can be expressed as a linear combination of~$\langle u_0\rangle_{\!,x}$ and~$\langle u_0\rangle_{\!,xx}$. Since~$\langle u_0\rangle$ is a generic mean wavefield propagating through a heterogeneous medium, however, $\langle u_0\rangle_{\!,x}$ and $\langle u_0\rangle_{\!,xx}$ are themselves linearly independent. As a result, their respective multipliers in~\eqref{Ptilde_Q} must vanish independently for each~$x$. By setting the multiplier therein of $\langle u_0\rangle_{\!,xx}$ to zero, we obtain the boundary value problem for $Q$ as
\begin{eqnarray} \label{CF_Q}
\begin{split}
& [G(P+Q_{,y})]_{,y} = \frac{\rho}{\varrho^{\mbox{\tiny{(0)}}}} \mu^{\mbox{\tiny{(0)}}}-G(1+P_{\!,y}),\quad y \in  Y_q \\*[-1mm]
& Q, ~ G(P+Q_{,y}) \quad Y\text{-periodic}; \quad \langle  Q \rangle = 0  \\*[-1mm]
& \llbracket Q\rrbracket =0, \quad \llbracket G(P+Q_{,y})\rrbracket =0,\quad   y \in \partial Y_q\backslash \partial Y
\end{split}
\end{eqnarray}
that holds for each~$x\in\mathbb{R}$. In this vein, the boundary value problem governing $\tilde{P}$ can be further identified as 
\begin{eqnarray} \label{CF_Ptilde}
\begin{split}
& [G(P_{,x}+\tilde{P}_{,y})]_{,y} = \frac{\rho}{\varrho^{\mbox{\tiny{(0)}}}} \mu^{\mbox{\tiny{(0)}}}_{,x}-[G(1+P_{\!,y})]_{,x}, \quad y \in  Y_q \\*[-1mm]
& \tilde{P}, ~ G(P_{,x}+\tilde{P}_{,y}) \quad Y\text{-periodic}; \quad \langle  \tilde{P} \rangle = 0  \\*[-1mm]
& \llbracket \tilde{P}\rrbracket =0, \quad \llbracket G(P_{,x}+\tilde{P}_{,y})\rrbracket =0, \quad   y \in \partial Y_q\backslash \partial Y. 
\end{split}
\end{eqnarray}

Proceeding with the analysis, we similarly find from~\eqref{WQR} the respective boundary value problems governing~$\tilde{Q}$, $R$ and $\tilde{R}$ to read
\begin{eqnarray} \label{CF_Qtilde}
\begin{split}
& [G(\tilde{P}+Q_{,x}+\tilde{Q}_{,y})]_{,y} =  \frac{\rho}{\varrho^{\mbox{\tiny{(0)}}}} (\eta+\mu^{\mbox{\tiny{(1)}}}_{,x}) - [G(P+Q_{,y})]_{,x}   
-G(P_{,x}+\tilde{P}_{,y}) \\*[-1mm] 
& \hspace*{55mm}+\rho \Big(P-\frac{\varrho^{\mbox{\tiny{(1)}}}}{\varrho^{\mbox{\tiny{(0)}}}}\Big)\Big(\frac{\mu^{\mbox{\tiny{(0)}}}_{,x}}{\varrho^{\mbox{\tiny{(0)}}}}+\big[\frac{\mu^{0}}{\varrho^{\mbox{\tiny{(0)}}}}\big]_{,x}\Big), \quad y \in  Y_q \\*[-1mm]
& \tilde{Q}, ~ G(\tilde{P}+Q_{,x}+\tilde{Q}_{,y}), \quad Y\text{-periodic}; \quad \langle  \tilde{Q} \rangle = 0  \\
& \llbracket \tilde{Q}\rrbracket =0, \quad \llbracket G(\tilde{P}+Q_{,x}+\tilde{Q}_{,y})\rrbracket =0, \quad   y \in \partial Y_q\backslash \partial Y,
\end{split}
\end{eqnarray}

\begin{eqnarray} \label{CF_R}
\begin{split}
& [G(Q+R_{,y})]_{,y}= \frac{\rho}{\varrho^{\mbox{\tiny{(0)}}}} \mu^{\mbox{\tiny{(1)}}} - G(P+Q_{,y})+ \rho 
\frac{\mu^{\mbox{\tiny{(0)}}}}{\varrho^{\mbox{\tiny{(0)}}}} \Big(P-\frac{\varrho^{\mbox{\tiny{(1)}}}}{\varrho^{\mbox{\tiny{(0)}}}}\Big), \quad y \in  Y_q \\*[-1mm]
& R, ~ G(Q+R_{,y}), \quad Y\text{-periodic}; \quad \langle  R \rangle = 0  \\*[-1mm]
& \llbracket R\rrbracket =0, \quad \llbracket G(Q+R_{,y})\rrbracket =0, \quad   y \in \partial Y_q\backslash \partial Y,
\end{split} 
\end{eqnarray}
and 
\begin{eqnarray} \label{CF_W}
\begin{split}
& [G(\tilde{R}_{,y}+\tilde{P}_{,x})]_{,y} = \frac{\rho}{\varrho^{\mbox{\tiny{(0)}}}} \eta _{,x} -[G(P_{,x}+\tilde{P}_{,y})]_{,x} + \rho \Big(P-\frac{\varrho^{\mbox{\tiny{(1)}}}}{\varrho^{\mbox{\tiny{(0)}}}}\Big)\big[\frac{\mu^{\mbox{\tiny{(0)}}}_{,x}}{\varrho^{\mbox{\tiny{(0)}}}}\big]_{,x}, \quad y \in  Y_q \\*[-1mm]
& \tilde{R}, \quad G(\tilde{R}_{,y}+\tilde{P}_{,x}), \quad Y\text{-periodic}; \quad \langle  \tilde{R} \rangle = 0  \\*[-1mm]
& \llbracket \tilde{R}\rrbracket =0, \quad \llbracket G(\tilde{R}_{,y}+\tilde{P}_{,x})\rrbracket =0, \quad   y \in \partial Y_q\backslash \partial Y,
\end{split}
\end{eqnarray}
for any given~$x\in\mathbb{R}$. 

\begin{remark}
Thanks to the results in~\citep{Bensoussan} (Chapter~2), one finds that the boundary value problems~\eqref{CF_P}--\eqref{CF_W} are well-posed, each featuring a field equation with the common principal part $[G(\cdot)_{,y}]_{,y}$. Therein, the dependence on the macroscopic variable~$x$ is injected via smooth coefficient variations~$G'(x)$ and~$\rho'(x)$ according to either~\eqref{HetMat11d} or~\eqref{HetMat21d}. As a result, the cell functions~$\{P,Q,R\}$ and~$\{\tilde{P},\tilde{Q},\tilde{R}\}$ are likewise smooth functions of~$x$, resulting in the smooth spatial variation of the effective coefficients~$\mu^{\mbox{\tiny{($k$)}}}, \varrho^{\mbox{\tiny{($k$)}}} (k=0,1,,\ldots)$, $\tilde{\varrho}^{\mbox{\tiny{(2)}}}$, $\eta,\phi$, and~$\psi$. In the context of numerical (e.g. finite element) implementation, such ``functionally-graded'' coefficients can then be sampled with suitable density in~$x$ and interpolated accordingly. In other words, even though the cell problems~\eqref{CF_P}--\eqref{CF_W} are dependent on the macroscopic variable~$x$, their solution would need to be sampled only over a relatively coarse grid.  
\end{remark} 

\subsection{Cell stresses and stress expansion} \label{Cell_Stress}

\noindent With the cell functions at hand, we next introduce the so-called cell stresses $\Sigma_j$ ($j=\overline{0,5}$) as 
\begin{eqnarray} \label{Cell_Sterss_Def}
\begin{split}
&   \Sigma_0(x,y) = \frac{G (1+P_{,y})}{\mu^{\mbox{\tiny{(0)}}}},   &   \Sigma_1(x,y) = \frac{G (P+Q_{,y})}{\mu^{\mbox{\tiny{(0)}}}}, \\
&   \Sigma_2(x,y) = \frac{G (P_{,x}+\tilde{P}_{,y})}{\mu^{\mbox{\tiny{(0)}}}},     &  \Sigma_3(x,y) = \frac{G (\tilde{R}_{,y}+\tilde{P}_{,x})}{\mu^{\mbox{\tiny{(0)}}}}, \\
&   \Sigma_4(x,y) = \frac{G (\tilde{P}+Q_{,x}+\tilde{Q}_{,y})}{\mu^{\mbox{\tiny{(0)}}}},    &   \Sigma_5(x,y) = \frac{G (Q+R_{,y})}{\mu^{\mbox{\tiny{(0)}}}},\\
\end{split}
\end{eqnarray}
By way of~\eqref{rho0_mu0} and~\eqref{CF_P}, one can show that $\Sigma_0(x,y)=1$; for clarity of discussion, however, we will retain this term ``as is'' wherever it appears. Thanks to~\eqref{sub}, \eqref{AsEx} and~\eqref{Cell_Sterss_Def}, the stress field~\eqref{stress1} affiliated with~$u(x,y)$ can be expanded as 
\begin{multline} \label{Stress_Expanded}
  \sigma(x,y) ~=~ \mu^{\mbox{\tiny{(0)}}} (x) \big[ \big(\Sigma_0(x,y)+\epsilon \hspace{0.5mm}\Sigma_2(x,y)+\epsilon^2 \hspace{0.5mm}\Sigma_3(x,y)\big)\langle u\rangle^{[2]}_{,x}   \\ 
 + \big(\epsilon \hspace{0.5mm}\Sigma_1(x,y)+\epsilon^2 \hspace{0.5mm}\Sigma_4(x,y)\big)\langle u\rangle^{[2]}_{,xx} + \epsilon^2 \hspace{0.5mm}\Sigma_5(x,y)\langle u\rangle^{[2]}_{,xxx} \big] + O(\epsilon^3).
\end{multline}

\section{Effective field equation in~\mbox{$\mathbb{R}^d$ $(d>1)$}} \label{bi_periodic} 

\noindent In this section, we apply the foregoing two-scale analysis to the original problem~\eqref{WaveEq} in~$\mathbb{R}^d$. For brevity, we show only the essential definitions and homogenization results. We start by introducing the ``fast'' spatial coordinate $\by=\eps^{-1}\bx$ and mappings 
\begin{eqnarray} \label{sub_2D}
\begin{split}
& u(\bx) \Rightarrow u(\bx,\by ),\qquad   G(\bx) \Rightarrow G(\bx,\by ),    \\ 
&\rho(\bx) \Rightarrow \, \rho(\bx,\by ), \qquad  \nabla  \Rightarrow \nabla _{\bx} +\frac{1}{\epsilon} \nabla _{\by }, 
\end{split}
\end{eqnarray}
which transforms~\eqref{WaveEq} into
\begin{eqnarray} \label{ReScEq_2D}
\epsilon^{-2}[\nabla_{\!\by } \!\cdot\! (G\nabla_{\!\by } u)] \,+\, \epsilon^{-1}[    \nabla_{\!\bx} \!\cdot\! ( G  \nabla_{\!\by } u ) + \nabla_{\!\by } \!\cdot\! ( G  \nabla_{\!\bx} u )] \,+\,  
[\rho \omega^2 u + \nabla_{\!\bx} \!\cdot\! ( G  \nabla_{\!\bx} u )] \;=\; 0.
\end{eqnarray}
We then pursue the ansatz 
\begin{eqnarray} \label{AsEx_2d}
u(\bx,\by)=\sum_{k=0}^{\infty} \epsilon^k u_k (\bx,\by) \quad \Longrightarrow \quad \langle u\rangle(\bx) = \int_{ Y} u(\bx,\by) \, \text{d}\by = \sum_{k=0}^{\infty} \epsilon^k \langle u_k\rangle (\bx). 
\end{eqnarray}
By retracing the steps of one-dimensional analysis, we specifically find that 
\begin{eqnarray} \label{u2_2d}
u_2 (\bx,\by) \;=\; \langle u_2\rangle  + \boldsymbol{P}(\bx,\by) \!\cdot\! \nabla  \langle u_1\rangle +\boldsymbol{Q}(\bx,\by) : \nabla  \nabla \langle u_0\rangle +\tilde{\boldsymbol{P}}(\bx,\by) \!\cdot\!  \nabla \langle u_0\rangle,  \quad \langle \boldsymbol{Q} \rangle = \langle \boldsymbol{P}\rangle = \langle \tilde{\boldsymbol{P}} \rangle =\bzero,
\end{eqnarray}
where~$\boldsymbol{P}, \boldsymbol{Q}$ and~$\tilde{\boldsymbol{P}}$ are \emph{tensorial} cell functions reading 
\[
\boldsymbol{P} = P_i \, \be_i, \qquad \boldsymbol{Q} = Q_{ij} \, \be_i \otimes\be_j, \qquad \tilde{\boldsymbol{P}} = \tilde{P}_i \,\be_i 
\]
in dyadic notation, which assumes implicit summation over repeated indexes $i,j = \overline{1,d}$. In~\eqref{u2_2d} and thereafter, we use symbol ``:'' to indicate $n$-tuple contraction between two $n$th order tensors ($n\geqslant 2$) producing a scalar. On further denoting by~$\bI_n$ the $n$th order symmetric identity tensor, the mean fields~$\langle u_k\rangle$ ($k=0,1,2$) featured in~\eqref{u2_2d} can be shown to satisfy the respective field equations
\begin{eqnarray} \label{O_1_Eq_2d}
\nabla\sip [\boldsymbol{\mu^{\mbox{\tiny{(0)}}}} (\bx)\sip\nabla \langle u_0\rangle ] \,+\, \varrho^{\mbox{\tiny{(0)}}}(\bx)\es \omega^2 \langle u_0\rangle \;=\; 0
\end{eqnarray}
where
\begin{eqnarray} \label{rho0_mu0_2d}
\boldsymbol{\mu^{\mbox{\tiny{(0)}}}}(\bx)=\langle    G ( \boldsymbol{I}_2  + \nabla_{\!\by } \boldsymbol{P} ) \rangle, \qquad 
\varrho^{\mbox{\tiny{(0)}}}(\bx) = \langle \rho \rangle;
\end{eqnarray}

\begin{multline} \label{O_eps_Eq_2d}
\nabla\sip[\boldsymbol{\mu^{\mbox{\tiny{(0)}}}} (\bx)\sip\nabla \langle u_1\rangle] \,+\, \varrho^{\mbox{\tiny{(0)}}}(\bx)\es \omega^2 \langle u_1\rangle \;=\;  - [\nabla\sip\boldsymbol{\eta}(\bx)+\boldsymbol{\varrho^{\mbox{\tiny{(1)}}}}(\bx)\es \omega^2] \!\cdot \! \nabla \langle u_0\rangle \\  - [\nabla \sip \boldsymbol{\mu^{\mbox{\tiny{(1)}}}}(\bx)+\boldsymbol{\eta }(\bx) ] \!:\! \nabla \nabla \langle u_0\rangle 
- \boldsymbol{\mu^{\mbox{\tiny{(1)}}}}(\bx) : \nabla \nabla \nabla \langle u_0\rangle,
\end{multline}
where
\begin{eqnarray} \label{xiC_2d}
\boldsymbol{\mu^{\mbox{\tiny{(1)}}}}(\bx) = \langle  G(\boldsymbol{I}_2 \otimes \boldsymbol{P} + \nabla_{\!\by } \boldsymbol{Q} ) \rangle, 
\qquad \boldsymbol{\varrho^{\mbox{\tiny{(1)}}}} (\bx)= \langle  \rho \boldsymbol{P} \rangle, \qquad 
\boldsymbol{\eta}(\bx)= \langle  G(\nabla_{\!\bx} \boldsymbol{P}+\nabla_{\!\by } \tilde{\boldsymbol{P}}) \rangle; 
\end{eqnarray}
and 
\begin{multline} \label{O_eps2_Eq_2d}
\nabla \sip [\boldsymbol{\mu^{\mbox{\tiny{(0)}}}}(\bx) \sip  \nabla \langle u_2\rangle] +\varrho^{\mbox{\tiny{(0)}}}(\bx)\es \omega^2 \langle u_2\rangle \;=\; - [\nabla \sip \boldsymbol{\eta } (\bx) + \boldsymbol{\varrho^{\mbox{\tiny{(1)}}}}(\bx)\es \omega^2] \sip \nabla \langle u_1\rangle  
- [\boldsymbol{\eta }(\bx) + \nabla \sip \boldsymbol{\mu^{\mbox{\tiny{(1)}}}} (\bx)] : \nabla \nabla \langle u_1\rangle \\  
- \boldsymbol{\mu^{\mbox{\tiny{(1)}}}}(\bx) : \nabla \nabla \nabla  \langle u_1\rangle   
- [\nabla\sip\boldsymbol{\phi}(\bx) + \tilde{\boldsymbol{\varrho}}^{\mbox{\tiny{(2)}}}(\bx)\es \omega^2] \sip \nabla \langle u_0\rangle  
 - [\nabla\sip\boldsymbol{\psi}(\bx)+\boldsymbol{\phi}(\bx) + \boldsymbol{\varrho^{\mbox{\tiny{(2)}}}}(\bx)\es\omega^2]:\nabla\nabla \langle u_0\rangle \\
- [ \boldsymbol{\psi}(\bx)+\nabla \sip \boldsymbol{\mu^{\mbox{\tiny{(2)}}}}(\bx)] \!:\! \nabla \nabla \nabla \langle u_0\rangle -  \boldsymbol{\mu^{\mbox{\tiny{(2)}}}} (\bx) : \nabla \nabla \nabla \nabla \langle u_0\rangle , 
\end{multline}
\noindent where
\begin{eqnarray} \label{phiC_2d}
\begin{split}
& \boldsymbol{\mu^{\mbox{\tiny{(2)}}}} (\bx) = \langle  G  (\boldsymbol{I}_2 \otimes  \boldsymbol{Q}+ \nabla_{\!\by } \boldsymbol{R} )\rangle, \qquad 
\boldsymbol{\varrho^{\mbox{\tiny{(2)}}}} (\bx) = \langle  \rho \boldsymbol{Q} \rangle, \qquad \tilde{\boldsymbol{\varrho}}^{\mbox{\tiny{(2)}}} (\bx)= \langle  \rho \tilde{\boldsymbol{P}} \rangle,  \\
&\boldsymbol{\phi} (\bx) = \langle  G(\nabla_{\!\bx} \tilde{\boldsymbol{P}}+\nabla_{\!\by } \tilde{\boldsymbol{R}}) \rangle, \qquad 
\boldsymbol{\psi} (\bx)= \langle  G(\nabla_{\!\bx}\boldsymbol{Q} +\boldsymbol{I}_2 \otimes  \tilde{\boldsymbol{P}} + \nabla_{\!\by } \tilde{\boldsymbol{Q}} ) \rangle.
\end{split}
\end{eqnarray}

\noindent From the weighted sum $\eqref{O_1_Eq_2d}+\eps\eqref{O_eps_Eq_2d}+\eps^2\eqref{O_eps2_Eq_2d}$, we obtain the effective field equation 
\begin{eqnarray} \label{MMEq_2d}
\begin{split}
& \nabla \sip  [\boldsymbol{\mu^{\mbox{\tiny{(0)}}}}(\bx) \sip  \nabla \langle u\rangle] + \varrho^{\mbox{\tiny{(0)}}}(\bx)\es \omega^2 \langle u\rangle \;=\;  \\
& \qquad -\epsilon \es \Big([\nabla \sip \boldsymbol{\eta}(\bx) + \boldsymbol{\varrho^{\mbox{\tiny{(1)}}}}(\bx) \omega^2]  \sip \nabla \langle u\rangle + [\boldsymbol{\eta}(\bx) + \nabla \sip  \boldsymbol{\mu^{\mbox{\tiny{(1)}}}}(\bx)] : \nabla \nabla \langle u\rangle + \boldsymbol{\mu^{\mbox{\tiny{(1)}}}}(\bx) : \nabla \nabla \nabla \langle u\rangle \Big)  \\
& \qquad - \epsilon^2 \big([\nabla \sip  \boldsymbol{\phi}(\bx) + \tilde{\boldsymbol{\varrho}}^{\mbox{\tiny{(2)}}}(\bx)\es \omega^2] \sip \nabla \langle u\rangle + [\nabla \sip \boldsymbol{\psi}(\bx) + \boldsymbol{\phi}(\bx) + \boldsymbol{\varrho^{\mbox{\tiny{(2)}}}}(\bx)\es \omega^2] : \nabla \nabla \langle u\rangle   \\
&\qquad\qquad\quad  + [\boldsymbol{\psi}(\bx) + \nabla \sip  \boldsymbol{\mu^{\mbox{\tiny{(2)}}}}(\bx)] : \nabla \nabla \nabla \langle u\rangle + \boldsymbol{\mu^{\mbox{\tiny{(2)}}}} : \nabla \nabla \nabla \nabla \langle u\rangle \Big) \,+\, O(\eps^3), 
\end{split}
\end{eqnarray}  
where the cell functions needed to compute the effective coefficients~$\boldsymbol{P}\in\mathbb{R}^d,\boldsymbol{Q}\in\mathbb{R}^{d\times d},\boldsymbol{R}\in\mathbb{R}^{d\times d\times d}, \tilde{\boldsymbol{P}}\in\mathbb{R}^d, \tilde{\boldsymbol{Q}}\in\mathbb{R}^{d\times d}$ and~$\tilde{\boldsymbol{R}}\in\mathbb{R}^{d}$ according to~\eqref{rho0_mu0_2d}, \eqref{xiC_2d} and~\eqref{phiC_2d} are given in~\ref{appen1}. 

\begin{remark} Taking into account the above tensorial character of the respective cell functions, one finds from~\eqref{rho0_mu0_2d}, 
\eqref{xiC_2d} and~\eqref{phiC_2d} that $\varrho^{\mbox{\tiny{(0)}}}$ is a scalar; $\boldsymbol{\rho^{\mbox{\tiny{(1)}}}}$ and $\tilde{\boldsymbol{\rho}}^{\mbox{\tiny{(2)}}}$ are vectors; $\boldsymbol{\mu^{\mbox{\tiny{(0)}}}}, \boldsymbol{\rho^{\mbox{\tiny{(2)}}}}, \boldsymbol{\eta }$ and  $\boldsymbol{\phi}$ are second-order tensors; $\boldsymbol{\mu^{\mbox{\tiny{(1)}}}}$ and~$\boldsymbol{\psi}$ are third-order tensors, and $\boldsymbol{\mu^{\mbox{\tiny{(2)}}}}$ is a fourth-order tensor. 
\end{remark}

\section{Effective boundary conditions for one-dimensional problems} \label{HBC}

\noindent The subject of effective boundary conditions is, even for periodic media, highly challenging problem that admits explicit formulation only under special circumstances, see e.g.~\citep{Cako2019} and references therein. One such amenable class of problems is one-dimensional wave motion in bounded periodic domains where the effective boundary conditions, when expanded up to the second order~\citep{Corna2019}, take the form of (elastic spring-like) \emph{Robin conditions}~\citep{Akin2005} written in terms of the mean field. This motivates our attempt to introduce the effective boundary conditions for \emph{quasi-periodic media} as described next. To our knowledge, this problem has not been considered in the literature. 

\begin{figure}[h!]  \vspace*{2mm}
\centering{\includegraphics[width=0.6\linewidth]{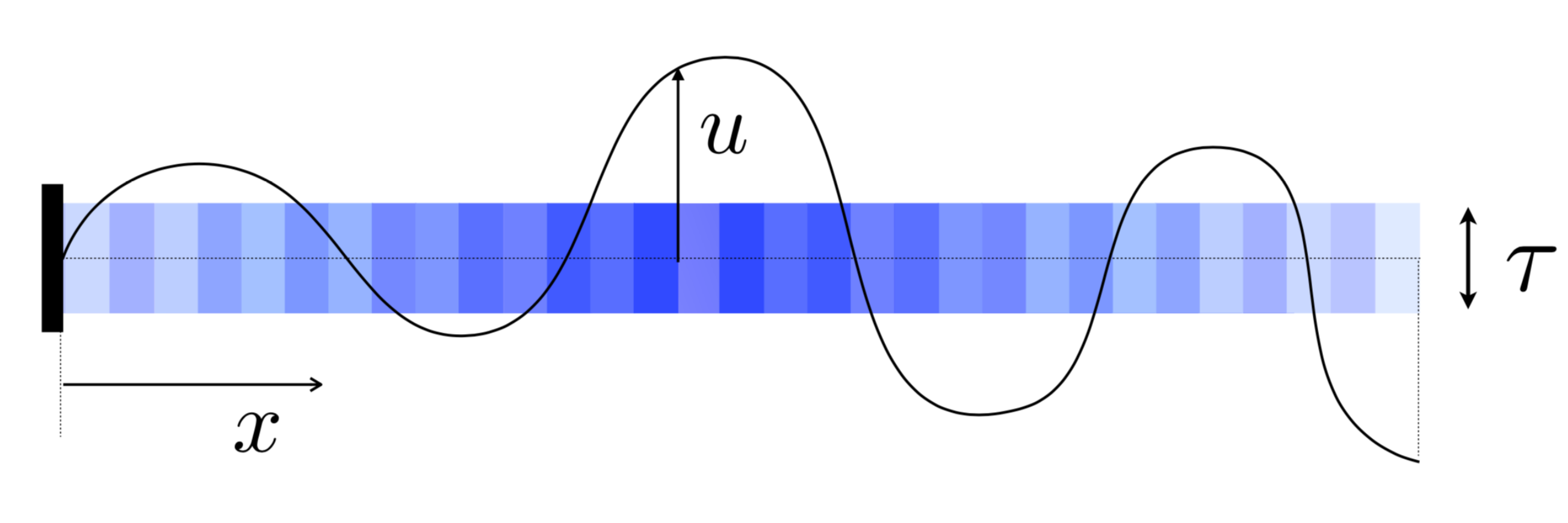}}\vspace*{-4mm}
\caption{Shear waves in a finite quasi-periodic medium.}
\label{1Dmedium}
\end{figure}

Let us consider a mixed boundary value problem (BVP) describing the shear wave motion in a quasi-periodic medium $Y=(0,1)$ that is fixed at~$x=0$ and subjected to time-harmonic shear traction $\tau$ at $x=1$, see Fig.~\ref{1Dmedium}. With reference to~\eqref{1D_Wave}, the BVP reads 

\begin{eqnarray} \label{PB}
\begin{aligned}
& \frac{\text{d} }{\text{d} x} \Big( G(x) \frac{\text{d} u}{\text{d} x} \Big) +\rho(x)\es \omega^2 u = 0, \qquad & x\in Y, \\
& u=0, &  x=0, \\
& \tau = G(x) \frac{\text{d} u}{\text{d} x}, & x=1,
\end{aligned}
\end{eqnarray}
where~$G$ and~$\rho$ are given by either~\eqref{HetMat1} or~\eqref{HetMat2}. 

\subsection{Effective field equation}

\noindent To facilitate the development of effective boundary conditions, it is next useful to rewrite the second-order effective field equation~(\ref{MMEq}) by expressing the featured third- and fourth-order derivatives in terms of their lower-order companions ($\langle u\rangle, \langle u\rangle_{\!,x}, \langle u\rangle_{\!,xx}$) using the $O(1)$ effective field equation~(\ref{O_1_Eq}) and the $O(\epsilon)$ effective field equation~(\ref{O_eps_Eq}). Accordingly, (\ref{MMEq}) can be rewritten compactly as 
\begin{alignat} {3}
\big\{E_5(x) + \omega^2 E_3(x) \big\} \langle u \rangle_{\!,xx} \,+\,
\big\{E_4(x) + \omega^2 E_2(x) \big\} \langle u \rangle_{\!,x} \,+\, \omega^2 E_1(x) \langle u\rangle \;=\; O(\eps^3), \qquad x \in \mathbb{R} \label{Reduced_Eq}
\end{alignat}
\noindent where
\begin{eqnarray} \label{EEE}
\begin{split}
& E_1 (x) = \varrho^{\mbox{\tiny{(0)}}}  +  \epsilon \bigg(\rho^{\mbox{\tiny{(1)}}} - \mu^{\mbox{\tiny{(1)}}} [\frac{\rho^{\mbox{\tiny{(0)}}}}{\mu^{\mbox{\tiny{(0)}}}} ]_{,x} \bigg)
  +  \epsilon^2 \bigg(    \mu^{\mbox{\tiny{(1)}}}  \Big\{  [ \frac{\rho^{\mbox{\tiny{(0)}}}}{\mu^{\mbox{\tiny{(0)}}}} ]_{,x}  \Big( \frac{\mu^{\mbox{\tiny{(1)}}}_{,x} + \eta^{\mbox{\tiny{(1)}}} }{\mu^{\mbox{\tiny{(0)}}}} + [\frac{\mu^{\mbox{\tiny{(1)}}} }{\mu^{\mbox{\tiny{(0)}}}} ]_{,x} \Big) + \\ & \hspace*{15mm} \frac{\mu^{\mbox{\tiny{(1)}}} }{\mu^{\mbox{\tiny{(0)}}}} \Big(   [ \frac{\rho^{\mbox{\tiny{(0)}}}}{\mu^{\mbox{\tiny{(0)}}}} ]_{,xx} -   [ \frac{\mu^{\mbox{\tiny{(0)}}}_{,x}}{\mu^{\mbox{\tiny{(0)}}}} ] [ \frac{\rho^{\mbox{\tiny{(0)}}}}{\mu^{\mbox{\tiny{(0)}}}} ]_{,x}\Big)  \Big\}    - [\frac{\varrho^{\mbox{\tiny{(0)}}}}{\mu^{\mbox{\tiny{(0)}}}}]_{,x}  \big(\psi+\mu^{\mbox{\tiny{(2)}}}_{,x}\big) + \mu^{\mbox{\tiny{(2)}}} \Big\{\frac{\mu^{\mbox{\tiny{(0)}}}_{,x}}{\mu^{\mbox{\tiny{(0)}}}} [\frac{\varrho^{\mbox{\tiny{(0)}}}}{\mu^{\mbox{\tiny{(0)}}}}]_{,x} -[\frac{\varrho^{\mbox{\tiny{(0)}}}}{\mu^{\mbox{\tiny{(0)}}}}]_{,xx} \Big\} \bigg)  \\ 
& E_2 (x) = -\epsilon \bigg( \mu^{\mbox{\tiny{(1)}}}   [ \frac{\rho^{\mbox{\tiny{(0)}}}}{\mu^{\mbox{\tiny{(0)}}}} ] \bigg)    +  \epsilon^2 \bigg(  \mu^{\mbox{\tiny{(1)}}} \Big \{ [ \frac{\rho^{\mbox{\tiny{(0)}}}}{\mu^{\mbox{\tiny{(0)}}}} ]  \Big( [\frac{\mu^{\mbox{\tiny{(1)}}}_{,x} + \eta^{\mbox{\tiny{(1)}}} }{\mu^{\mbox{\tiny{(0)}}}} ] + [\frac{\mu^{\mbox{\tiny{(1)}}} }{\mu^{\mbox{\tiny{(0)}}}} ]_{,x} \Big) + [\frac{\mu^{\mbox{\tiny{(1)}}} }{\mu^{\mbox{\tiny{(0)}}}}]  \Big( 2[\frac{\rho^{\mbox{\tiny{(0)}}}}{\mu^{\mbox{\tiny{(0)}}}} ]_{,x} - [ \frac{\mu^{\mbox{\tiny{(0)}}}_{,x}}{\mu^{\mbox{\tiny{(0)}}}} ][ \frac{\rho^{\mbox{\tiny{(0)}}}}{\mu^{\mbox{\tiny{(0)}}}} ] \Big) \\ & \hspace*{15mm}  - [\frac{\rho^{\mbox{\tiny{(1)}}} }{\mu^{\mbox{\tiny{(0)}}}} ]_{,x} \Big \} +   \tilde{\varrho}^{\mbox{\tiny{(2)}}}  -\frac{\varrho^{\mbox{\tiny{(0)}}}}{\mu^{\mbox{\tiny{(0)}}}} \big(\psi +\mu^{\mbox{\tiny{(2)}}}_{,x}\big) +\mu^{\mbox{\tiny{(2)}}} \Big\{ (\frac{\mu^{\mbox{\tiny{(0)}}}_{,x}}{\mu^{\mbox{\tiny{(0)}}}})(\frac{\varrho^{\mbox{\tiny{(0)}}}}{\mu^{\mbox{\tiny{(0)}}}}) - 2 [\frac{\varrho^{\mbox{\tiny{(0)}}}}{\mu^{\mbox{\tiny{(0)}}}}]_{,x} \Big\} \bigg )  \\
& E_3 (x) =  \epsilon^2 \bigg(\mu^{\mbox{\tiny{(1)}}} \Big \{ [\frac{\rho^{\mbox{\tiny{(0)}}}}{\mu^{\mbox{\tiny{(0)}}}} ] [\frac{\mu^{\mbox{\tiny{(1)}}} }{\mu^{\mbox{\tiny{(0)}}}}] - [\frac{\rho^{\mbox{\tiny{(1)}}} }{\mu^{\mbox{\tiny{(0)}}}} ] \Big \} + \varrho^{\mbox{\tiny{(2)}}} - \mu^{\mbox{\tiny{(2)}}} (\frac{\varrho^{\mbox{\tiny{(0)}}}}{\mu^{\mbox{\tiny{(0)}}}}) \bigg)    \\
& E_4 (x) =  \mu^{\mbox{\tiny{(0)}}}_{,x} +   \epsilon \bigg( \eta^{\mbox{\tiny{(1)}}}_{,x} - \mu^{\mbox{\tiny{(1)}}} [\frac{\mu^{\mbox{\tiny{(0)}}}_{,x}}{\mu^{\mbox{\tiny{(0)}}}}]_{,x}   \bigg)   +  \epsilon^2 \bigg(  \phi_{,x} -   [\frac{\mu^{\mbox{\tiny{(0)}}}_{,x}}{\mu^{\mbox{\tiny{(0)}}}}]_{,x} \big(  \psi + \mu^{\mbox{\tiny{(2)}}}_{,x} \big) + \mu^{\mbox{\tiny{(2)}}} \Big\{ (\frac{\mu^{\mbox{\tiny{(0)}}}_{,x}}{\mu^{\mbox{\tiny{(0)}}}}) [\frac{\mu^{\mbox{\tiny{(0)}}}_{,x}}{\mu^{\mbox{\tiny{(0)}}}}]_{,x} - [\frac{\mu^{\mbox{\tiny{(0)}}}_{,x}}{\mu^{\mbox{\tiny{(0)}}}}]_{,xx} \Big\}  \\ & \hspace*{15mm}  + \mu^{\mbox{\tiny{(1)}}} \Big \{ [\frac{\mu^{\mbox{\tiny{(0)}}}_{,x}}{\mu^{\mbox{\tiny{(0)}}}}]_{,x} \Big(  [ \frac{\mu^{\mbox{\tiny{(1)}}}_{,x} + \eta^{\mbox{\tiny{(1)}}} }{\mu^{\mbox{\tiny{(0)}}}} ]   + [ \frac{\mu^{\mbox{\tiny{(1)}}} }{\mu^{\mbox{\tiny{(0)}}}} ]_{,x}  \Big) -  [\frac{\eta^{\mbox{\tiny{(1)}}}_{,x}}{\mu^{\mbox{\tiny{(0)}}}} ]_{,x} + [\frac{\mu^{\mbox{\tiny{(1)}}} }{\mu^{\mbox{\tiny{(0)}}}}] \Big( [ \frac{\mu^{\mbox{\tiny{(0)}}}_{,x}}{\mu^{\mbox{\tiny{(0)}}}} ]_{,xx} - [ \frac{\mu^{\mbox{\tiny{(0)}}}_{,x}}{\mu^{\mbox{\tiny{(0)}}}} ] [ \frac{\mu^{\mbox{\tiny{(0)}}}_{,x}}{\mu^{\mbox{\tiny{(0)}}}} ]_{,x} \Big)   \Big\}     \bigg)    \\
& E_5 (x) =  \mu^{\mbox{\tiny{(0)}}} \,+\,     \epsilon \bigg(   \eta^{\mbox{\tiny{(1)}}} +\mu^{\mbox{\tiny{(1)}}}_{,x}  -\mu^{\mbox{\tiny{(1)}}} [ \frac{\mu^{\mbox{\tiny{(0)}}}_{,x}}{\mu^{\mbox{\tiny{(0)}}}} ]  \bigg )    \,+\,   \epsilon^2 \bigg(  \psi_{,x} + \phi -\frac{\mu^{\mbox{\tiny{(0)}}}_{,x}}{\mu^{\mbox{\tiny{(0)}}}} \big(  \psi + \mu^{\mbox{\tiny{(2)}}}_{,x} \big) + \mu^{\mbox{\tiny{(2)}}} \Big\{ [\frac{\mu^{\mbox{\tiny{(0)}}}_{,x}}{\mu^{\mbox{\tiny{(0)}}}}]^2 -  2[\frac{\mu^{\mbox{\tiny{(0)}}}_{,x}}{\mu^{\mbox{\tiny{(0)}}}}]_{,x} \Big\}   \\ & \hspace*{15mm} +  \mu^{\mbox{\tiny{(1)}}}  \Big \{[\frac{\mu^{\mbox{\tiny{(0)}}}_{,x}}{\mu^{\mbox{\tiny{(0)}}}} ]  \Big(  [ \frac{\mu^{\mbox{\tiny{(1)}}}_{,x} + \eta^{\mbox{\tiny{(1)}}} }{\mu^{\mbox{\tiny{(0)}}}} ] + [ \frac{\mu^{\mbox{\tiny{(1)}}} }{\mu^{\mbox{\tiny{(0)}}}} ]_{,x} \Big)  - [\frac{\eta^{\mbox{\tiny{(1)}}}_{,x}}{\mu^{\mbox{\tiny{(0)}}}} ] - [\frac{\mu^{\mbox{\tiny{(1)}}}_{,x} + \eta^{\mbox{\tiny{(1)}}} }{\mu^{\mbox{\tiny{(0)}}}} ]_{,x} + [ \frac{\mu^{\mbox{\tiny{(1)}}} }{\mu^{\mbox{\tiny{(0)}}}} ] \Big( 2 [ \frac{\mu^{\mbox{\tiny{(0)}}}_{,x}}{\mu^{\mbox{\tiny{(0)}}}} ]_{,x} - [ \frac{\mu^{\mbox{\tiny{(0)}}}_{,x}}{\mu^{\mbox{\tiny{(0)}}}} ]^2 \Big)  \Big \}      \bigg)
\end{split}
\end{eqnarray}

\subsection{Zeroth-order model} \label{HBC_0thModel}

\noindent By virtue of~\eqref{resEx}, \eqref{u0} and~\eqref{Stress_Expanded}, one finds that the zeroth-order, \emph{single-scale} approximations of the displacement field~$u(x)$ and stress field~$\sigma(x)$ solving~\eqref{PB} can be written respectively as 
\begin{eqnarray} \label{PB_u0}
\begin{split}
& u^{\mbox{\tiny{[0]}}}(x) = \langle u\rangle^{\mbox{\tiny{[0]}}}, \\
&  \sigma^{\mbox{\tiny{[0]}}} (x) = \mu^{\mbox{\tiny{(0)}}} (x) \, \langle u\rangle^{\mbox{\tiny{[0]}}}_{,x} \; \Sigma_0(x,x/\epsilon), 
\end{split}
\end{eqnarray}
where $\Sigma_0(x,x/\epsilon)=1$ as examined earlier. By requiring the \emph{mean-field} equation~\eqref{Reduced_Eq} and the \emph{local} boundary conditions in~\eqref{PB} to be satisfied up to~$O(1)$, the mean field $\langle u\rangle^{\mbox{\tiny{[0]}}}$ can be shown to satisfy the effective BVP
\begin{eqnarray} \label{BVP_u0}
\begin{aligned}
& \bar{E}_{5}(x) \langle u\rangle^{\mbox{\tiny{[0]}}}_{,xx} \,+\, \bar{E}_{4}(x) \langle u\rangle^{\mbox{\tiny{[0]}}}_{,x} + \omega^2 \bar{E}_{1}(x) \langle u\rangle^{\mbox{\tiny{[0]}}} \:=\: 0, \qquad & x\in Y \\
& \langle u\rangle^{\mbox{\tiny{[0]}}} =0, & x=0, \\
&   \langle u\rangle^{\mbox{\tiny{[0]}}}_{,x}  = \frac{\tau}{\mu^{\mbox{\tiny{(0)}}}(1)}, & x=1,
\end{aligned}
\end{eqnarray}
where $\bar{E}_{1}$, $\bar{E}_{4}$, and $\bar{E}_{5}$ denote respectively the truncations of $E_{1}$, $E_{4}$, and $E_{5}$ that retain terms up to $O(1)$.

\subsection{First-order model} \label{HBC_1stModel}

\noindent By pursuing the analysis similar to that in Section \ref{HBC_0thModel}, we next proceed with the first-order homogenized model. On recalling~\eqref{resEx} and making use of~\eqref{u1} and~\eqref{Cell_Sterss_Def}, the first-order approximations of $u(x)$ and $\sigma(x)$ satisfying~\eqref{PB} can be written as     
\begin{eqnarray} \label{PB_u1}
\begin{split}
& u^{\mbox{\tiny{[1]}}}(x) = \langle u\rangle^{\mbox{\tiny{[1]}}} + \epsilon\hspace{0.5mm} P(x,x/\epsilon) \langle u\rangle^{\mbox{\tiny{[1]}}}_{,x}, \\
&  \sigma^{\mbox{\tiny{[1]}}}(x) = \mu^{\mbox{\tiny{(0)}}} (x) \big[ (\Sigma_0(x,x/\epsilon)+\epsilon \hspace{0.5mm}\Sigma_2(x,x/\epsilon))\langle u\rangle^{\mbox{\tiny{[1]}}}_{,x} +  \epsilon \hspace{0.5mm} \Sigma_1(x,x/\epsilon) \langle u\rangle^{\mbox{\tiny{[1]}}}_{,xx} \big] ,\\
\end{split}
\end{eqnarray}

\noindent Then, by requiring the \emph{mean-field} equation~\eqref{Reduced_Eq} and the \emph{local} boundary conditions in~\eqref{PB} to be satisfied up to~$O(\eps)$, the mean field $\langle u\rangle^{\mbox{\tiny{[1]}}}$ can be shown to satisfy the effective BVP 
\begin{eqnarray} \label{BVP_u1}
\begin{aligned}
& \hat{E}_{5}(x) \langle u \rangle^{\mbox{\tiny{[1]}}}_{,xx} \,+\, \hat{E}_{4}(x) \langle u \rangle^{\mbox{\tiny{[1]}}}_{,x}  +  \omega^2 \hat{E}_{1}(x) \langle u \rangle^{\mbox{\tiny{[1]}}} \:=\: 0, \qquad  &  x\in Y, \\
& \langle u\rangle^{\mbox{\tiny{[1]}}} + \epsilon \hspace{0.5mm} P(0,0) \langle u \rangle^{\mbox{\tiny{[1]}}}_{,x}=0, & x=0, \\
&  \Big\{-\epsilon \hspace{0.5mm}\Sigma_1(1,1)\hspace{0.5mm} \omega^2 \frac{\bar{E}_1(1)}{\bar{E}_5(1)}\Big\} \langle u\rangle^{\mbox{\tiny{[1]}}}  +  \Big\{\Sigma_0(1,1) + \epsilon \hspace{0.5mm} (\Sigma_2(1,1) - \Sigma_1(1,1) \frac{\bar{E}_4(1)}{\bar{E}_5(1)}) \Big\} \langle u\rangle^{\mbox{\tiny{[1]}}}_{,x}  \:=\: \frac{\tau}{\mu^{\mbox{\tiny{(0)}}}(1)},  &  x=1,
\end{aligned}
\end{eqnarray}
where $\hat{E}_{j}$ $(j=\overline{1,5})$ denotes the truncation of $E_{j}$ that discards $O(\eps^2)$ correction terms. Note that in~\eqref{BVP_u1}, we utilized the featured field equation $\langle u \rangle_{\!,xx}^{\mbox{\tiny{[1]}}}= -\hat{E}_5^{-1}(\omega^2 \hat{E}_1 \langle u \rangle^{\mbox{\tiny{[1]}}} + \hat{E}_4 \langle u \rangle_{\!,x}^{\mbox{\tiny{[1]}}})$ in the stress boundary condition in order to obtain Robin-type boundary condition at each end.

\subsection{Second-order model} \label{HBC_2ndModel}

\noindent Proceeding with the second-order homogenization, one can similarly write the second-order approximations of $u(x)$ and $\sigma(x)$ solving~\eqref{PB} as
\begin{eqnarray} \label{PB_u2}
\begin{split}
& u^{\mbox{\tiny{[2]}}}(x) = \langle u\rangle^{\mbox{\tiny{[2]}}} + [ \epsilon\hspace{0.5mm} P(x,x/\epsilon)  +  \epsilon^2\hspace{0.5mm} \tilde{P}(x,x/\epsilon) ] \langle u\rangle^{\mbox{\tiny{[2]}}}_{,x}+  \epsilon^2\hspace{0.5mm} Q(x,x/\epsilon) \langle u\rangle^{\mbox{\tiny{[2]}}}_{,xx} , \\
&  \sigma^{\mbox{\tiny{[2]}}}(x) = \mu^{\mbox{\tiny{(0)}}} (x) \big[ \big(\Sigma_0(x,x/\epsilon)+\epsilon \hspace{0.5mm}\Sigma_2(x,x/\epsilon)+\epsilon^2 \hspace{0.5mm}\Sigma_3(x,x/\epsilon)\big)\langle u\rangle^{\mbox{\tiny{[2]}}}_{,x}  \\ 
& \hspace*{14mm} + \big(\epsilon \hspace{0.5mm}\Sigma_1(x,x/\epsilon)+\epsilon^2 \hspace{0.5mm}\Sigma_4(x,x/\epsilon)\big)\langle u\rangle^{\mbox{\tiny{[2]}}}_{,xx} + \epsilon^2 \hspace{0.5mm}\Sigma_5(x,x/\epsilon)\langle u\rangle^{\mbox{\tiny{[2]}}}_{,xxx} \big] ,\\
\end{split}
\end{eqnarray}
where the second-order mean field $\langle u\rangle^{\mbox{\tiny{[2]}}}$ solves the  BVP
\begin{eqnarray} \label{BVP_u2}
\begin{aligned}
& \big\{E_{5}(x)  +  \omega^2 E_{3}(x)\big\} \langle u \rangle^{\mbox{\tiny{[2]}}}_{,xx} \,+\, \big\{E_{4}(x)   +  \omega^2 E_{2}(x)\big\} \langle u \rangle^{\mbox{\tiny{[2]}}}_{,x}  +  \omega^2 E_{1}(x) \langle u \rangle^{\mbox{\tiny{[2]}}} \:=\: 0, &  x\in Y, \\
& \Big\{ 1-\epsilon^2 Q(0,0) \mathcal{D}(0) \mathcal{F}(0) \Big\} \langle u\rangle^{\mbox{\tiny{[2]}}} + \Big\{\epsilon \hspace{0.5mm} P(0,0) + \epsilon^2 \tilde{P}(0,0) - \epsilon^2 Q(0,0) \mathcal{D}(0) \mathcal{E}(0)  \Big\}  \langle u \rangle^{\mbox{\tiny{[2]}}}_{,x}=0, & x=0, \\
& \textcolor{black}{\Big\{\mathcal{C}(1,1)\mathcal{H}_1(1) - \mathcal{D}(1) \mathcal{B}(1,1) \mathcal{F}(1) \Big\} \langle u\rangle^{\mbox{\tiny{[2]}}}  +  \Big\{ \mathcal{A}(1,1) +\mathcal{C}(1,1)\mathcal{H}_2(1) - \mathcal{D}(1) \mathcal{B}(1,1)  \mathcal{E}(1) \Big\} \langle u\rangle^{\mbox{\tiny{[2]}}}_{,x}  = \frac{\tau}{\mu^{\mbox{\tiny{(0)}}}(1)},} &  x=1, \\
\end{aligned}
\end{eqnarray}
where
\begin{eqnarray} \label{BVP_u2_Coeff}
\begin{split}
& \mathcal{A}(x,x/\epsilon) = \Sigma_0(x,x/\epsilon)+\epsilon \,\Sigma_2(x,x/\epsilon)+\epsilon^2 \,\Sigma_3(x,x/\epsilon), \quad
\mathcal{B} (x,x/\epsilon)=  \epsilon \,\Sigma_1(x,x/\epsilon)+\epsilon^2 \,\Sigma_4(x,x/\epsilon), \\ & \mathcal{C} (x,x/\epsilon)= \epsilon^2 \,\Sigma_5(x,x/\epsilon), \quad
 \mathcal{D} (x)= \big(\omega^2 E_{3}(x) +E_{5}(x)\big)^{-1}, \quad 
\mathcal{E} (x)=  \omega^2 E_{2}(x) + E_{4}(x), \quad
\mathcal{F} (x)= \omega^2 E_{1}(x) ,\\
& \mathcal{H}_1 (x)= -\mathcal{D}_{,x} \mathcal{F}+\mathcal{D}^2 \mathcal{E} \mathcal{F}  -\mathcal{D} \mathcal{F}_{,x},\\
& \mathcal{H}_2 (x)= -\mathcal{D}_{,x} \mathcal{E} -\mathcal{D} \mathcal{E}_{,x} + \mathcal{D}^2 \mathcal{E}^2  - \mathcal{D} \mathcal{F}(x),\\
\end{split}
\end{eqnarray}

\noindent Note that in~\eqref{BVP_u2}, we utilized the featured field equation $\langle u \rangle_{\!,xx}^{\mbox{\tiny{[2]}}}= -\mathcal{D} \big[\mathcal{F} \langle u \rangle^{\mbox{\tiny{[2]}}} + \mathcal{E} \langle u \rangle_{\!,x}^{\mbox{\tiny{[2]}}} \big]$ and its derivative $\langle u \rangle^{\mbox{\tiny{[2]}}}_{,xxx} = \mathcal{H}_1\langle u \rangle^{\mbox{\tiny{[2]}}} + \mathcal{H}_2 \langle u \rangle^{\mbox{\tiny{[2]}}}_{,x}$ in the stress boundary condition in order to obtain Robin-type boundary condition at each end.  

\begin{remark}
As can be seen from~\eqref{BVP_u2}, a second-order mean-field approximation of the BVP~\eqref{PB} entails (i)~second-order field equation with smooth coefficients, and (ii) Robin-type boundary conditions. The foregoing analysis can be easily generated to situations where the domain terminates within some~$\eps Y$, which then affects the constant coefficients specifying the boundary conditions. 
\end{remark}

\section{Numerical results} \label{BL}

\noindent To illustrate the utility of the foregoing homogenization framework, we consider one-dimensional wave motion~(\ref{1D_Wave}) in a quasi-periodic medium endowed with bilaminate microstructure. For generality, we consider both (a) waves in an \emph{unbounded} domain~$\mathbb{R}$ -- in terms of the heterogeneity-induced wave dispersion, and (b)~waves in a \emph{bounded} domain~$Y$ -- in terms of the waveforms generated by prescribed boundary excitation. With reference to~\eqref{HetMat1}, we consider several examples or quasi-periodic structures endowed with macroscopic variation 
\begin{alignat}{3}
\text{linear variation:} \quad  & G'(x) = 1 + \gamma_{\mbox{\tiny G}} \hh\hh x, \qquad 
&& \rho'(x) = 1 + \gamma_{\rho} \hh\hh x  \label{grlin}  \\ 
\text{sine~variation:} \quad  & G'(x) = 1 + \gamma_{\mbox{\tiny G}}\hh \sin(2\pi x + \beta_{\mbox{\tiny G}}), \quad  && \rho'(x) = 1 + \gamma_{\rho}\hh \sin(2\pi x + \beta_{\rho}) \label{grsin} 
\end{alignat}
and \emph{piecewise-constant} microscopic fluctuation 
\begin{equation}\label{grmicro}
G''(y) \,=\, \delta_{\mbox{\tiny G}} \big(\!-\!1 + 2H(y-\alpha) \big), \qquad \rho''(y) \,=\, \delta_{\rho} \big(\!-\!1 + 2H(y-\alpha) \big), \qquad 0<y<1
\end{equation}
where~$H$ is the Heaviside function and $\alpha\in(0,1)$, $\delta_{\mbox{\tiny G}}$, and $\delta_{\rho}$ are prescribed constants. 

\begin{remark}
As stated earlier, all quantities in this study are assumed to be normalized with respect to ``suitable dimensional basis''. In this section, the latter is given by the triplet $\{G_{\circ},~ \rho_{\circ},~L_{\circ}\}$, where~$G_{\circ}$ and~$\rho_{\circ}$ correspond to the respective constant terms in the macroscopic variations~\eqref{grlin}--\eqref{grsin} of the shear modulus and mass density, while~$L_\circ$ signifies the physical length of the unit cell~$Y$. 
\end{remark}

\subsection{Effective coefficients} \label{Coeff_Homo}

\begin{figure}[h!]  \vspace*{-1mm}
\centering{\includegraphics[width=0.95\linewidth]{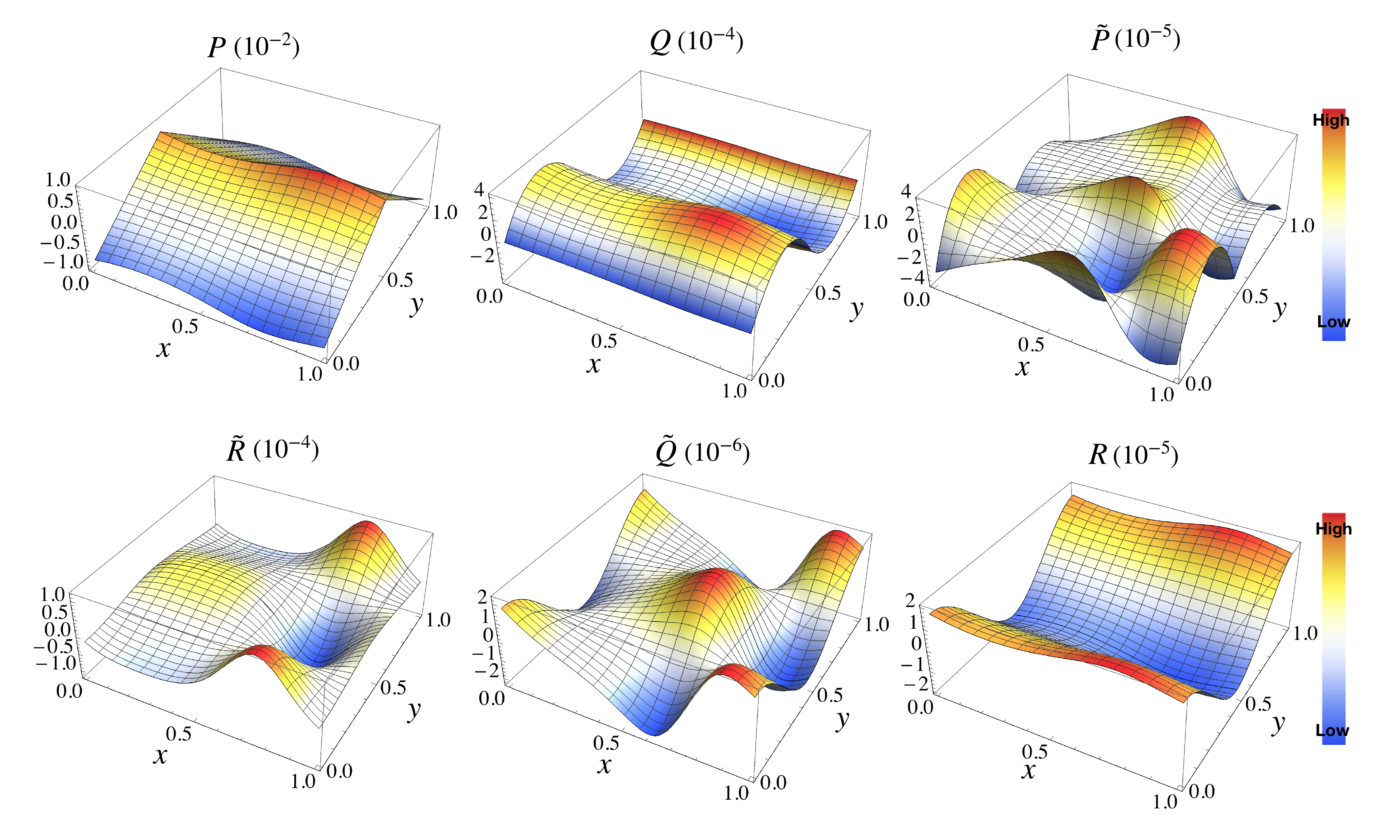}}\vspace*{-5mm}
\caption{Cell functions $P$, $Q$, $\tilde{P}$, $\tilde{R}$, $\tilde{Q}$ and~$R$: \SRC{Material~3}, sinusoidal profile.} 
\label{Cell_Funcx}
\end{figure}

\begin{figure}[h!]  \vspace*{-1mm}
\centering{\includegraphics[width=0.95\linewidth]{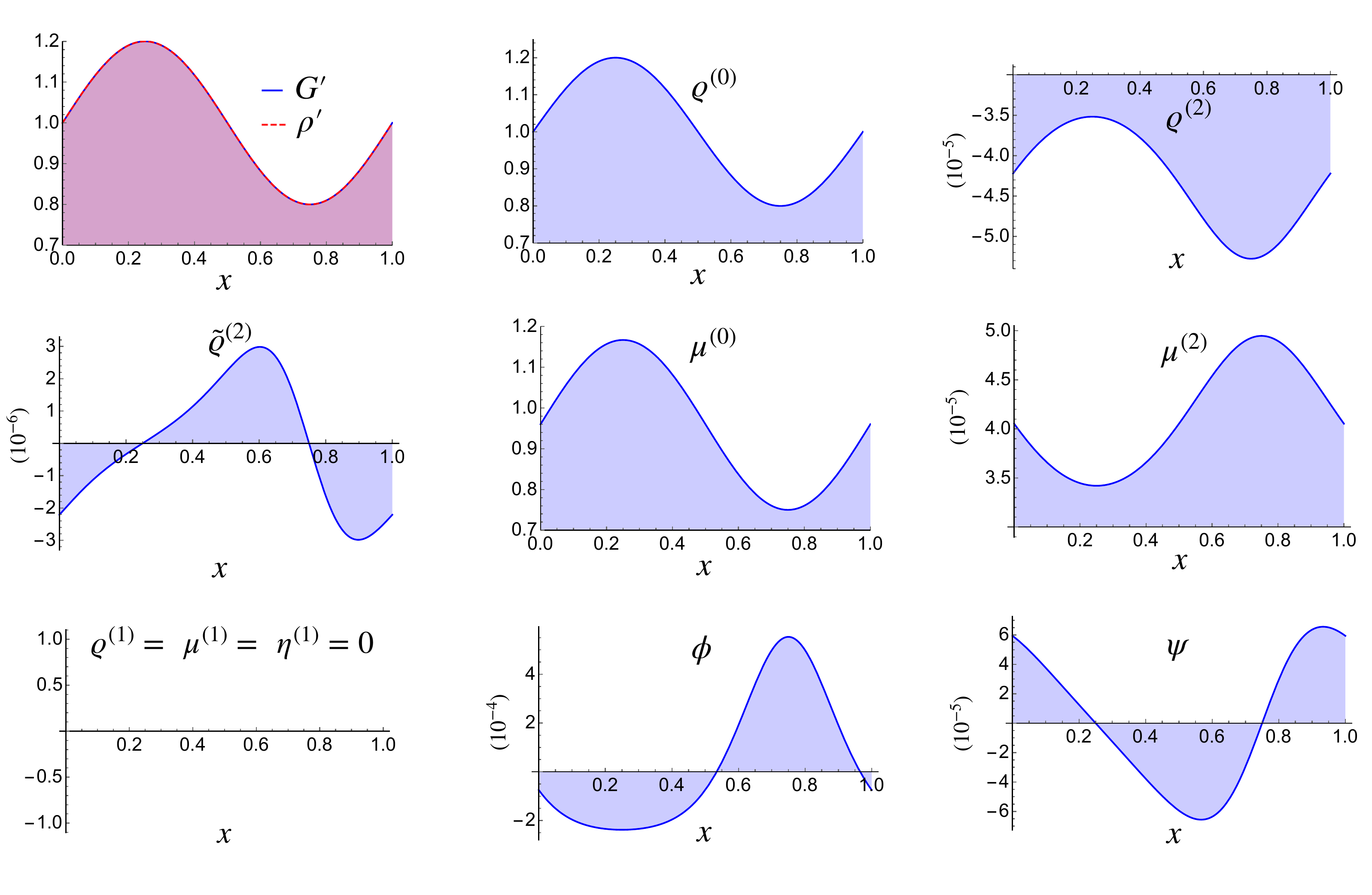}}\vspace*{-7mm}
\caption{Macroscopic medium variations~$(G',\rho')$ and effective coefficients $\varrho^{\mbox{\tiny{(0)}}}$, $\varrho^{\mbox{\tiny{(2)}}}$,   $\tilde{\varrho}^{\mbox{\tiny{(2)}}}$, $\mu^{\mbox{\tiny{(0)}}}$,  $\mu^{\mbox{\tiny{(2)}}}$, $\varrho^{\mbox{\tiny{(1)}}}$, $\mu^{\mbox{\tiny{(1)}}}$, $\eta$, $\phi$, and~$\psi$ featured by the second-order model~\eqref{MMEq_2d}: Material~3, sinusoidal profile.} 
	\label{Effective_Coeff}
\end{figure}

\noindent By applying the analysis from Section \ref{MuSc} to the class~\eqref{grlin}--\eqref{grmicro} of quasi-periodic media and solving the germane boundary value problems (\ref{CF_P})--(\ref{CF_W}), the effective coefficients of homogenization can be computed from~(\ref{rho0_mu0}), (\ref{xiC}), and~(\ref{phiC}). An analytical solution in terms of the cell functions $P,Q,\tilde{P},\tilde{R},\tilde{Q}$ and $R$ is sought via the symbolic manipulation platform Mathematica. Note that for \emph{any macroscopic variation} of the quasi-periodic medium~\eqref{HetMat1} with bilaminate microstructure~\eqref{grmicro}, one can explicitly compute $\varrho^{\mbox{\tiny{(0)}}}, \mu^{\mbox{\tiny{(0)}}}, \varrho^{\mbox{\tiny{(1)}}}$ and $\mu^{\mbox{\tiny{(1)}}}$ as
\begin{eqnarray} \label{effec_coeff}
\begin{split}
&\varrho^{\mbox{\tiny{(0)}}}(x) = \rho'(x) + \delta_{\mbox{\tiny G}} (1-2\alpha_{\mbox{\tiny G}}), \qquad 
\mu^{\mbox{\tiny{(0)}}}(x) = \frac{(G'(x)-\delta_{\mbox{\tiny G}}) (G'(x)+\delta_{\mbox{\tiny G}})}{\alpha(G'(x)+\delta_{\mbox{\tiny G}}) + (1-\alpha)(G'(x)-\delta_{\mbox{\tiny G}})},  \\
&\varrho^{\mbox{\tiny{(1)}}}(x) = 0, \qquad \mu^{\mbox{\tiny{(1)}}}(x) = 0,
\end{split}
\end{eqnarray}   
see also~\citep{Fish2001,Wautier} in the context of periodic media. Expressions for the remaining effective coefficients such as~$\eta(x)$ and~$\phi(x)$ in the second-order model~\eqref{MMEq}, however, are rather lengthy and will not be reported. Instead, a Mathematica code for their evaluation (assuming~\eqref{grlin}--\eqref{grmicro}) is provided as electronic supplementary material. From~\eqref{effec_coeff} it is also interesting to observe that $\mu^{\mbox{\tiny{(1)}}}$ and $\varrho^{\mbox{\tiny{(1)}}}$ vanish identically, which is a well known property of all periodic structures~\citep{Fish2001, Wautier}. With reference to Table~\ref{Table1} summarizing the example material profiles considered,  Fig.~\ref{Cell_Funcx} plots the featured cell functions for \SRC{Material~3}, while Fig.~\ref{Effective_Coeff} shows the distribution of the corresponding effective coefficients. 

\begin{table}[h!]
\caption{\label{Table1} Quasi-periodic profiles used in numerical simulations.} 
\begin{center}
\begin{tabular}{|l|c|c|c|c|c|c|c|c|} \hline
{$G'(x)$ and~$\rho'(x)$}  & $\gamma_G$ &$\delta_G$ & $\beta_G$ & $\gamma_{\rho} $  & $\delta_{\rho} $ & $\beta_{\rho}$ &$\alpha$\\ \hline\hline
Material $1$ :  \eqref{grsin}  & 1/5  & 3/5 & 0  & 1/5   &1/25  & 0 &  1/2\\ \hline
Material $2$ :  \eqref{grsin}  & 1/5  & 1/5 & 0  &0      & 0      & 0 &  1/2\\ \hline
Material $3$ :  \eqref{grsin}  & 1/5  & 1/5 &  0 & 1/5  & 1/5   & 0 & 1/2\\ \hline			
Material $4$ :  \eqref{grsin}  & 1/5 &  2/5 & $\pi/2$ & 1/5 & 1/5 & $\pi/2$ & 1/2\\ \hline
Material $5$ :  \eqref{grlin}  & $2 \pi/5$ & 2/5  & -- &$2 \pi/5$ & 1/5  & -- & 1/2\\ \hline
\end{tabular}
\end{center} 
\end{table}

As reported in~\citep{Wautier,Guillot}, higher-order effective coefficients for both periodic and quasi-periodic structures are typically small in magnitude relative to their leading-order companions, and Fig.~\ref{Effective_Coeff} confirms this observation. For the material profiles considered, it will be demonstrated that such behavior results in only a modest correction of the \emph{phase} of the asymptotic solution, but potentially significant \emph{amplitude} corrections.

\subsection{Wave dispersion} \label{dispersion}

\begin{figure}[h!]  \vspace*{-1mm}
\centering{\includegraphics[width=0.95\linewidth]{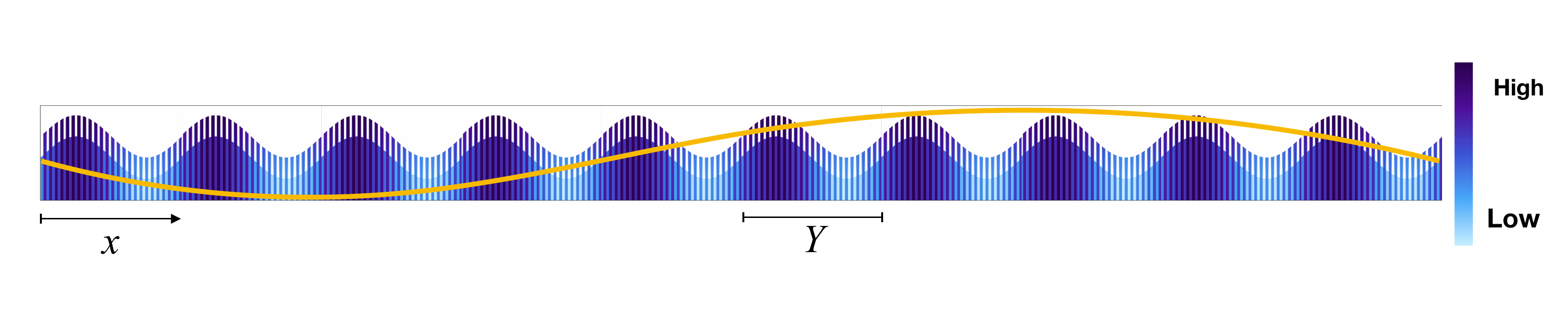}}\vspace*{-9mm}
\caption{Long-wavelength shear waves propagating through an infinite periodic medium endowed with sinusoidal macroscopic profile~\eqref{grsin} and bilaminate microstructure~\eqref{grmicro}.}
\label{Mega_Cell}
\end{figure}

\noindent We next consider the native wave equation~\eqref{WaveEq} for de facto \emph{periodic media} where (i) $G'$ and~$\rho'$ are $Y\!$-periodic, and (ii)~$\eps = n^{-1}$ ($n\in\mathbb{Z}^+$). In this setting, we pursue the dispersion analysis via the Floquet-Bloch approach~\citep{Floq1883} by seeking a solution in the form 
\begin{eqnarray} \label{Bloch}
u(x) \:=\: \tilde{u}(x)\hh e^{i k x}, \qquad \tilde{u}:  Y\text{-periodic}
\end{eqnarray}
where $Y\!=\!(0,1)$ as before, and~$k$ is the wavenumber. Letting~$\eps = n^{-1}$ for some~$n\!>\!1$, it is clear that $G=G'+G''$ and~$\rho=\rho'+\rho''$ are also~$Y\!$-periodic. In this setting, we are in position to \emph{homogenize} the \emph{``macrocell''}~$Y$ (containing~$n$ periods of the microstructural variation) and compare the dispersion relationship computed in this way with numerical simulations of the native Floquet-Bloch problem given by~\eqref{WaveEq} and~\eqref{Bloch}. In what follows, we present the dispersion results for a sinusoidal macroscopic profile~\eqref{grsin}, illustrated schematically in Fig.~\ref{Mega_Cell}.

By way of~\eqref{Bloch}, governing equation~(\ref{1D_Wave}) in~$\mathbb{R}$ reduces to the macrocell problem 
\begin{eqnarray} \label{BlochSol}
\begin{split}
& \frac{\text{d}}{\text{d}x_{\!k}}   \Big(G(x) \frac{\text{d}\tilde{u}}{\text{d}x_{\!k}} \Big) + \rho(x) \omega^2  \tilde{u} = 0, \quad \quad x \in Y \\
& \tilde{u}|_{x=0} = \tilde{u}|_{x=1}, \qquad G  \frac{\text{d}\tilde{u}}{\text{d} x_{\!k}} |_{x=0} = -    G  \frac{\text{d}\tilde{u}}{\text{d}x_{\!k}} |_{x=1}, 
\end{split}
\end{eqnarray}
where $\text{d}/\text{d}x_{\!k} := \text{d}/\text{d}x + i k$. Hereon, we refer to a numerical solution of the boundary value problem~(\ref{BlochSol}) as the ``exact'' solution. 

From~\eqref{Bloch}, it follows that 
\begin{eqnarray} \label{Bloch_Mean}
\langle u \rangle (x) \;=\; \langle \tilde{u} \rangle (x) \hh e^{i k x}.
\end{eqnarray}
Assuming the material profile~\eqref{grsin}--\eqref{grmicro} with $\eps = n^{-1}$ ($n>1$), one can make use of~\eqref{Bloch_Mean} together with the continuity and $Y$-periodicity of $E_3$ and $E_5$, to obtain the restriction of~\eqref{Reduced_Eq} to~$Y$ as
\begin{eqnarray} \label{Reduced_Eq_B}
\begin{split}
\omega^2 \big\{ E_6 (x) \langle \tilde{u} \rangle + E_7 (x) & \langle \tilde{u} \rangle_{\!,x} + E_3 (x) \langle \tilde{u} \rangle_{\!,xx} \big\} \;+\;  E_8 (x) \langle \tilde{u} \rangle + E_9 (x) \langle \tilde{u} \rangle_{\!,x} + E_5 (x) \langle \tilde{u} \rangle_{\!,xx} =0, \quad x \in Y \\
& \langle\tilde{u}\rangle|_{x=0} = \langle\tilde{u}\rangle|_{x=1}, \qquad  \langle \tilde{u} \rangle_{\!,x}|_{x = 0} = \langle \tilde{u} \rangle_{\!,x} |_{x = 1},
\end{split}
\end{eqnarray}
where
\begin{eqnarray} \label{EEE_Bloch}
\begin{split}
& E_6 (x) = E_1 (x)+ (ik) E_2(x) + (ik)^2 E_3(x),   \\
& E_7 (x) = E_2 (x)+ 2(ik) E_3(x), \\
& E_8 (x) =  (ik) E_4(x) + (ik)^2 E_5(x), \\
& E_9 (x) = E_4 (x)+ 2(ik) E_5(x). 
\end{split}  
\end{eqnarray}
Since the macroscopic variation of the medium is now $Y$-periodic by design, we exploit the framework of Section~3 to obtain its 
(constant) leading-order effective coefficients, and we use the affiliated (linear) dispersion relationship as a baseline in the ensuing simulations.

\begin{figure}[h!]  \vspace*{0mm}
\centering{\includegraphics[width=0.44\linewidth]{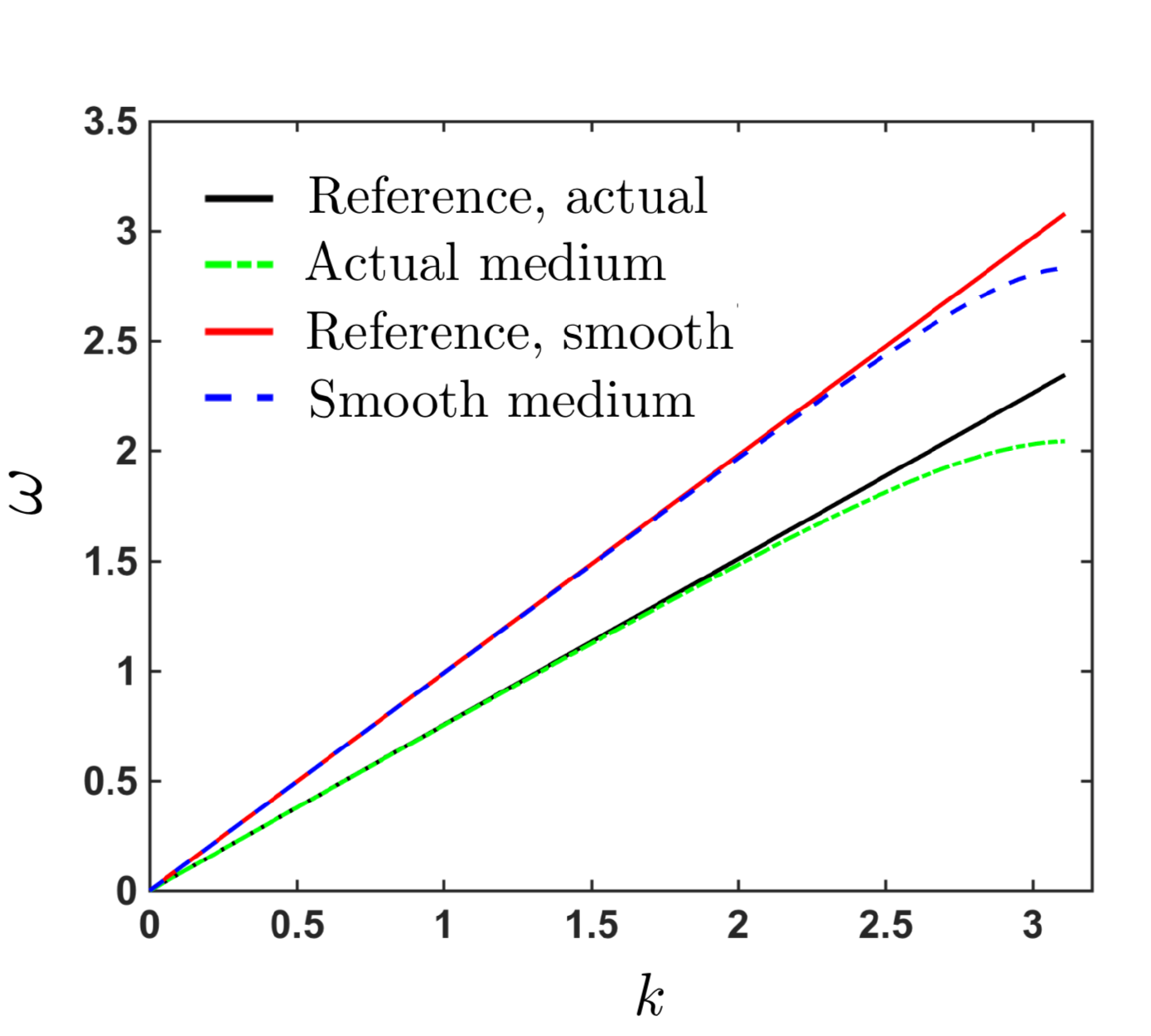}}\vspace*{-3mm}
\caption{Exact dispersion relationship due to~(\ref{BlochSol}) for an infinite periodic medium composed of Material~1 ($G=G'+G''$ and~$\rho=\rho'+\rho''$, green solid line) versus that for its \emph{microstructure-free} companion ($G=G'$ and~$\rho=\rho'$, blue dashed line). The fully homogenized, non-dispersive descriptions of~$Y$ for both media (black solid line for Material 1 and red solid line for its microstructure-free counterpart) are included as baselines.}
	\label{SmootVSac}
\end{figure}

To illustrate the analysis, we consider the wave dispersion in a ``sinusoidal'' medium~\eqref{grsin} endowed with bilaminate microstructure~\eqref{grmicro} and $\eps=1/50$. To examine the performance of the effective model, the dispersion curves for the ``exact'' solution and its homogenized counterparts are obtained via the Floquet-Bloch approach applied to the ``macrocell'' $Y$. As a point of reference, Fig.~\ref{SmootVSac} shows the ``exact'' dispersion curve in the first pass band (acoustic branch) for a medium composed of \SRC{Material~1}. To help parse the  effects of macroscopic  and microscopic heterogeneities on the overall behavior, we also include the dispersion curve for a \emph{microstructure-free} medium with~$G=G'$ and~$\rho=\rho'$. From the display, one observes that the microscopic medium fluctuations have a remarkable effect on the overall wave dispersion.

\begin{figure}[h!]  \vspace*{2mm}
\centering{\includegraphics[width=0.46\linewidth]{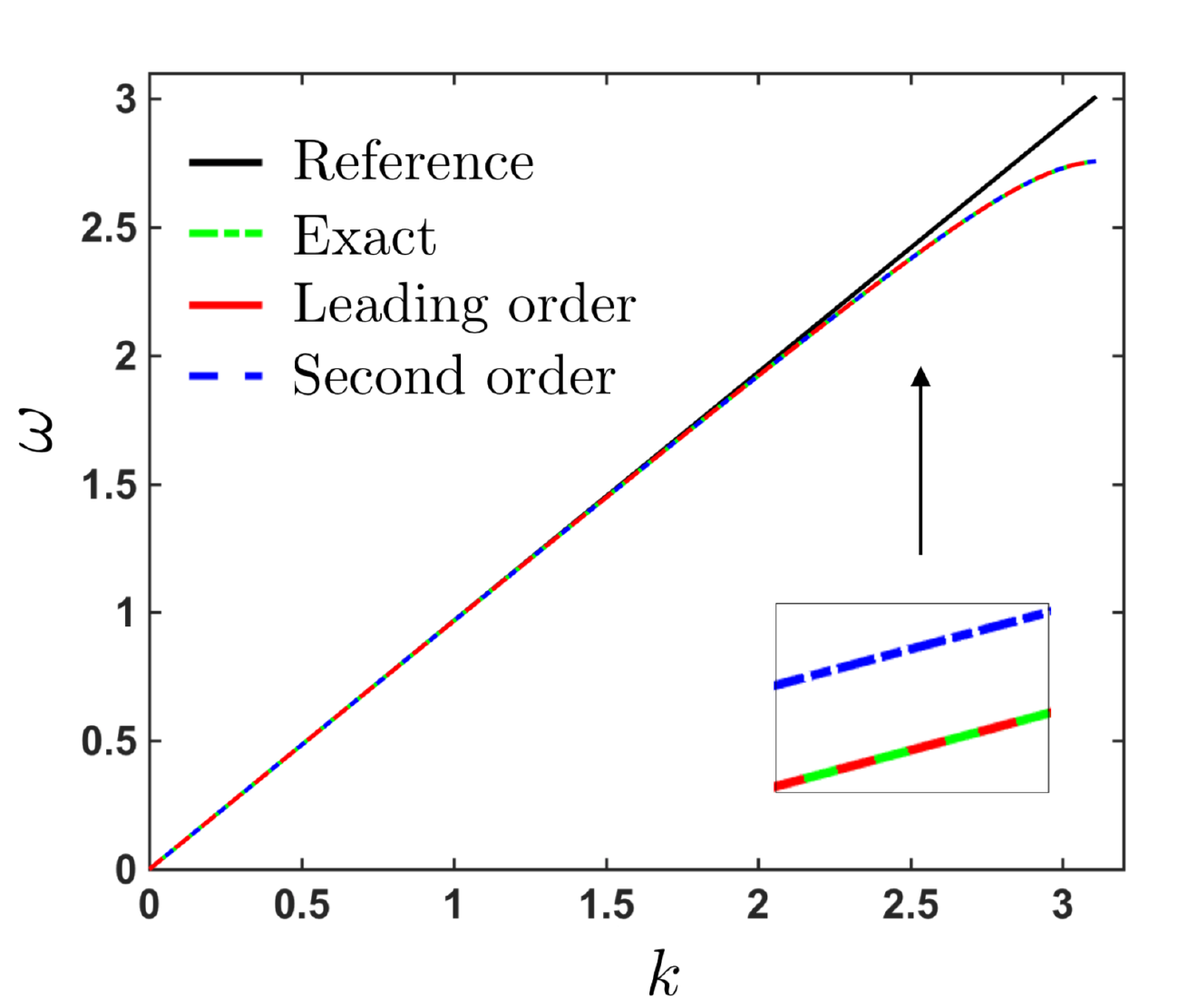}}\vspace*{-4mm}
\caption{Exact dispersion relationship due to~\eqref{BlochSol} for an infinite periodic medium composed of Material~3 (green solid line) versus its zeroth-order approximation (dot-dashed line) and second-order approximation (dashed line) according to~\eqref{Reduced_Eq_B}.}
	\label{Real_disp}
\end{figure}

\begin{figure}[h!]  \vspace*{4mm}
	\centering{\includegraphics[width=0.9\linewidth]{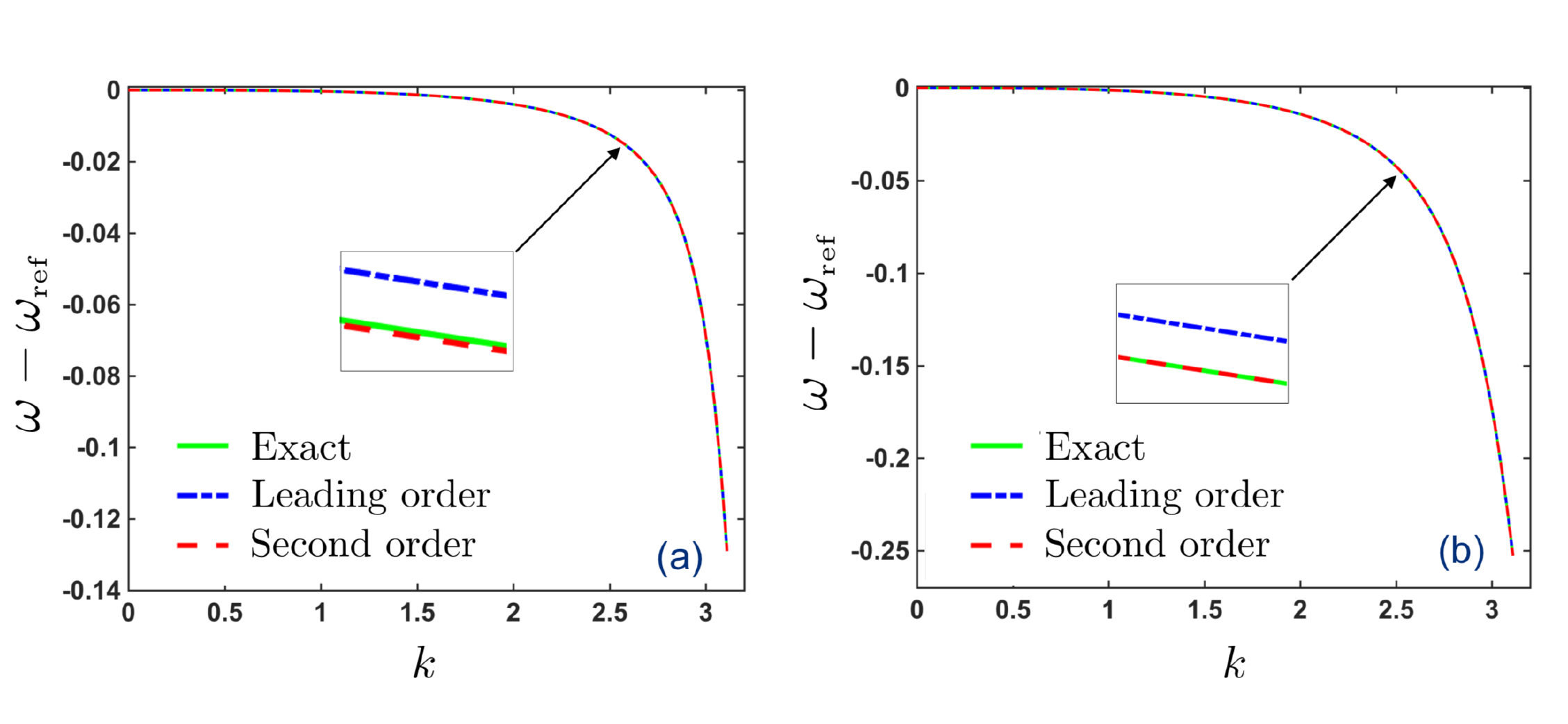}}\vspace*{-3mm}
	\caption{Separation from the reference non-dispersive model $\omega_{\mbox{\tiny{ref}}}(k)$: exact solution (green dashed line), zeroth-order approximation (blue dashed-dot line), and second-order approximation (red dashed line) for an infinite periodic medium composed of (a) \SRC{Material 2}, and (b) \SRC{Material 3}.}
	\label{modified_disp}
\end{figure}

Fig.~\ref{Real_disp} shows the dispersion curves in the first pass band for \SRC{Material~3}, and compares the exact model with its zeroth- and second-order asymptotic approximations stemming from~\eqref{Reduced_Eq_B}. As a point of reference, included in the diagram is a linear dispersion relationship for the leading-order (non-dispersive) homogenized macrocell~$Y$. The results suggest that the homogenized model is capable of capturing the exact solution with high accuracy. To highlight the dispersion effects, in Fig.~\ref{modified_disp} we recast the dispersion relationship in the terms of \emph{separation} from the reference (non-dispersive) relationship \SRC{for Material~2 and Material~3}. From the displays, one observes that even the leading-order, homogenized model of the quasi-periodic medium stemming from~\eqref{Reduced_Eq_B} captures the exact dispersion with high accuracy, which puts in question the utility of the higher-order corrections. As will be seen shortly, however, the conclusion changes drastically when considering the solution of a BVP in terms of actual waveforms -- that contain both amplitude and phase information.

\begin{figure}[h!]  \vspace*{8mm}
	\centering{\includegraphics[width=0.9\linewidth]{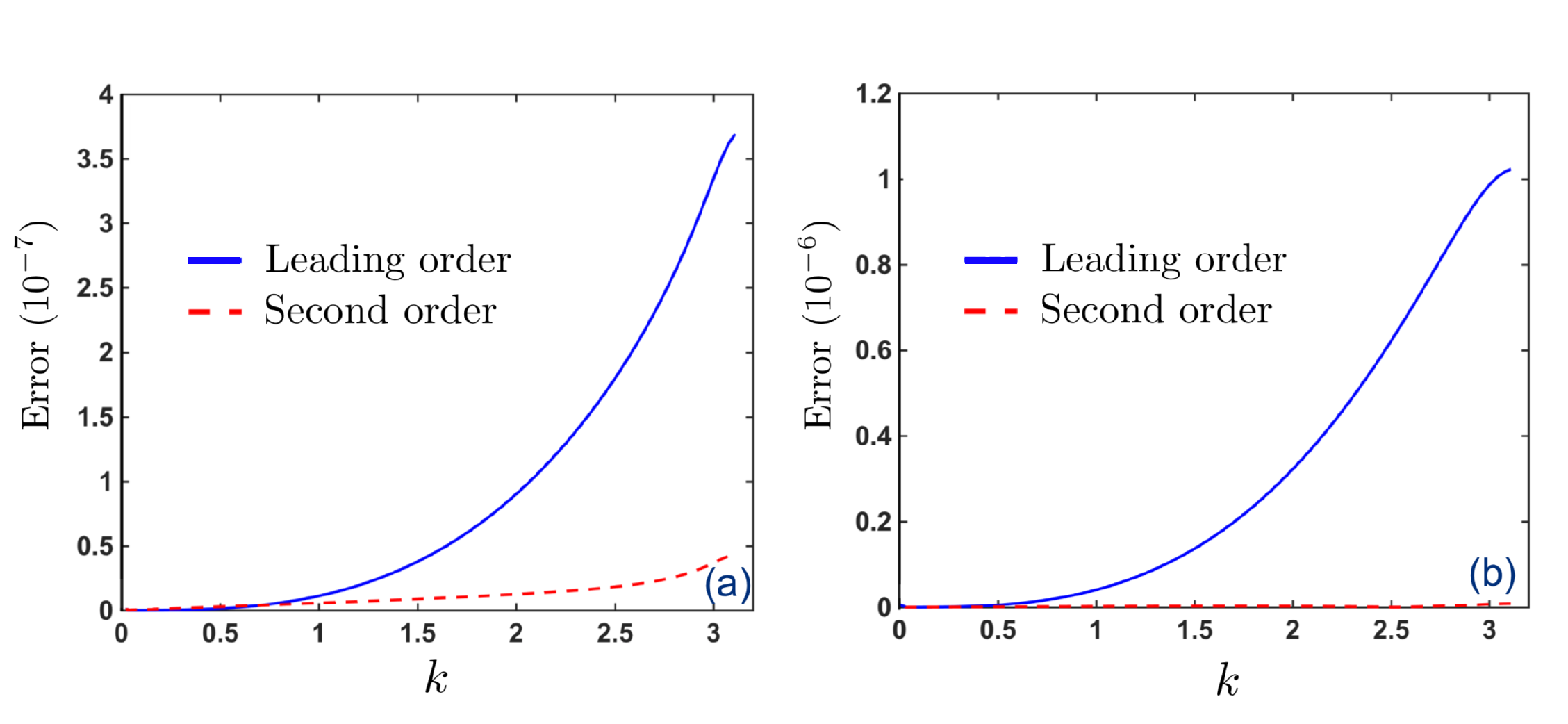}}\vspace*{-4mm}
	\caption{Relative error~\eqref{Error} of the asymptotic approximation ($j=0,2$) for an infinite periodic medium composed of (a)~\SRC{Material~2}, and (b)~\SRC{Material~3}.}
	\label{Error_WRT_Exact}
\end{figure}

For a complete insight into dispersive characteristics of the effective model, one can study the relative approximation error
\begin{eqnarray} \label{Error}
\text{Error} \;=\; \frac{|\omega^{(j)} - \omega^{(e)}|}{\|\omega^{(e)}\|_{2}}, \qquad j = 0,2
\end{eqnarray}
where~$\omega^{(e)}$ and~$\omega^{(j)}$ denote respectively the exact and $j$th-order homogenized solution, while $\|\cdot\|_{2}$ signifies the~$L^2$-norm computed over the (positive half of the) first Brillouin zone. In this setting, Fig.~\ref{Error_WRT_Exact} compares the relative error of the zeroth- and second-order effective models in describing the exact dispersion relationship for \SRC{Material~2 and Material~3}, which brings ``under microscope'' the improved fidelity brought about by the asymptotic correction. 

\subsection{Boundary value problem} \label{HBVP}

\noindent To complete the study, we consider an effective solution of the BVP~\eqref{PB} examined in Section~\ref{HBC} for the class of quasi-periodic media given by~\eqref{grlin}--\eqref{grmicro}. An ``exact'' solution of this problem is evaluated numerically via the propagator matrix approach using a high density of homogeneous sub-lamina to mimic the spatial variation of \SRC{$G=G'+G''$} and $\rho=\rho'+\rho$ (see~\citep{Wautier} for details in the context of periodic media). With reference to Table~\ref{Table1} specifying the example material profiles, Table~\ref{Table2} completes the list of input parameters required to simulate the BVP. 
\begin{table}[h!]
\caption{\label{Table2} Parameters used for numerical simulation of the BVP~\eqref{PB} and its homogenized approximations~\eqref{BVP_u0}, \eqref{BVP_u1}, and~\eqref{BVP_u2}.} 
\begin{center}
\begin{tabular}{|c|c|c|c|} \hline
Example  & Material & $\omega$ &$\epsilon \hspace{1mm}$ \\ \hline\hline
Ex$1$    &  $4$ & $ \pi ^2$   & 1/20 \\ \hline
Ex$2$    &  $4$ & $ 3 \pi ^2$ & 1/40  \\ \hline
Ex$3$    &  $5$ & $ 2 \pi^2 $ & 1/20 \\ \hline
\end{tabular}
\end{center} 
\end{table}

\begin{remark}
With reference to Fig.~\ref{Real_disp} describing the wave dispersion in an infinite periodic medium~$\mathbb{R}$ composed of Material~3 (and thus that of Material~4), excitation frequencies listed in Table~\ref{Table2} may appear as being ``too high'' in that they are located beyond the first pass band. In this regard, it is important to recognize that the problems examined in Section~\ref{dispersion} and Section~\ref{HBVP} are fundamentally different. Specifically, Section~\ref{dispersion} considers a macroscopically-periodic medium, that lends itself to the concept of the ``macroscopic'' Brillouin zone and allows for the computation of dispersion diagrams such as that in Fig.~\ref{Real_disp}. In contrast, this section is concerned with quasi-periodic media of finite extent, that are incompatible with the Bloch-wave representation~\eqref{Bloch} and the notion of the Brillouin zone. Instead, the key limitation on~$\omega$ in our general study of quasi-periodic media (see Section~\ref{OandA}) is that the lengthscale of microscopic fluctuations, $\eps Y$, is much smaller than the apparent wavelength. As will be seen shortly, all ensuing simulations meet this criterion by a safe margin. Concerning the remaining restriction on~$\omega$ from Section~\ref{OandA}, namely that the excitation frequency resides inside the first apparent ``pass band'', our numerical simulations show that increasing~$\omega$ beyond the values listed in Table~\ref{Table2} may lead to the creation of an apparent band gap, manifested by an exponential decay of wave amplitude away from the loaded end, $x=1$. Similarly, a reduction in~$\omega$ relative to the values listed in Table~\ref{Table2} can be shown (numerically) to result in a diminished error between that ``exact'' solution and its asymptotic approximation. 
\end{remark}

\begin{figure}[h!]  \vspace*{-1mm}
\centering{\includegraphics[width=0.82\linewidth]{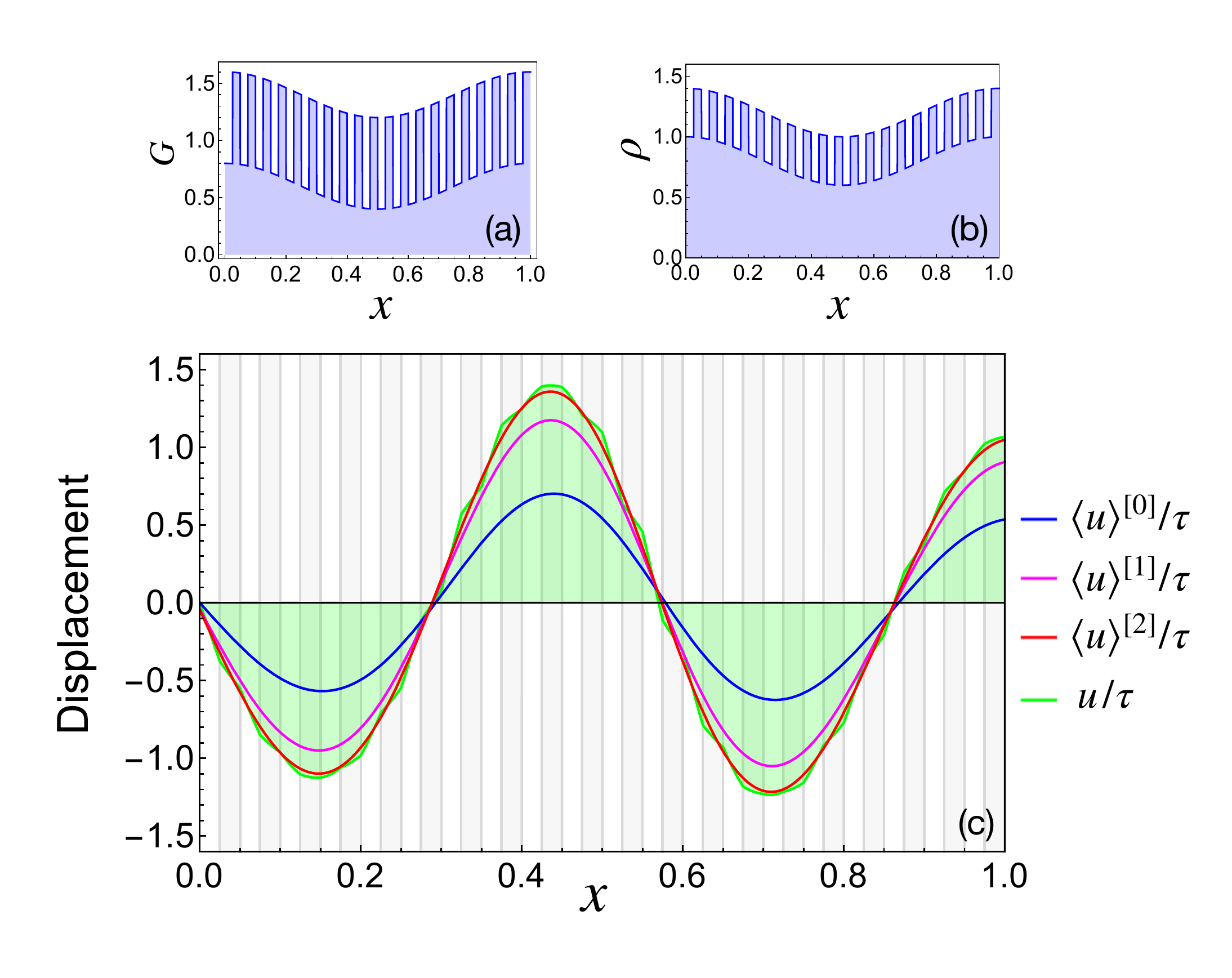}}\vspace*{-10mm}
\caption{Example Ex1: (a) Shear modulus profile, (b) mass density profile, and (c) ``exact'' wave motion $u$ versus homogenized mean fields $\langle u \rangle^{\mbox{\tiny{[0]}}}$, $\langle u \rangle^{\mbox{\tiny{[1]}}}$, and $\langle u \rangle^{\mbox{\tiny{[2]}}}$ for the BVP~\eqref{PB}.}
\label{PB_EX1_U}
\end{figure}

\begin{figure}[h!]  \vspace*{-1mm}
\centering{\includegraphics[width=0.78\linewidth]{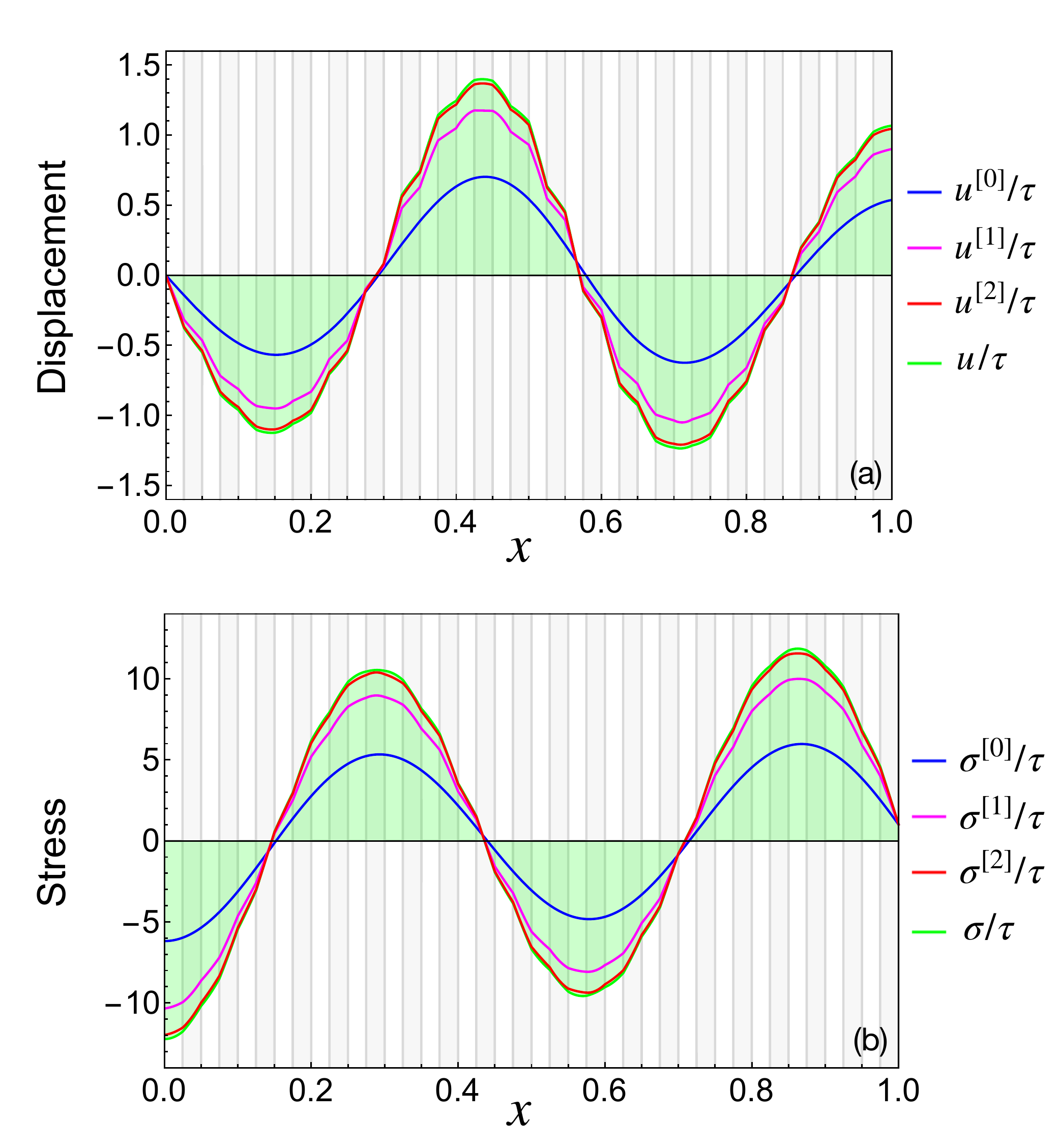}}\vspace*{-5mm}
\caption{Example Ex1: (a) ``exact'' wave motion $u$ versus homogenized approximations $u^{\mbox{\tiny{[0]}}}$, $u^{\mbox{\tiny{[1]}}}$, and $u^{\mbox{\tiny{[2]}}}$; (b) Exact stress filed $\sigma(x)$ versus homogenized approximations $\sigma^{\mbox{\tiny{[0]}}}$, $\sigma^{\mbox{\tiny{[1]}}}$, and $\sigma^{\mbox{\tiny{[2]}}}$, for the BVP~\eqref{PB}.}
\label{PB_EX1_us}
\end{figure}

\begin{figure}[h!]  \vspace*{-1mm}
\centering{\includegraphics[width=0.84\linewidth]{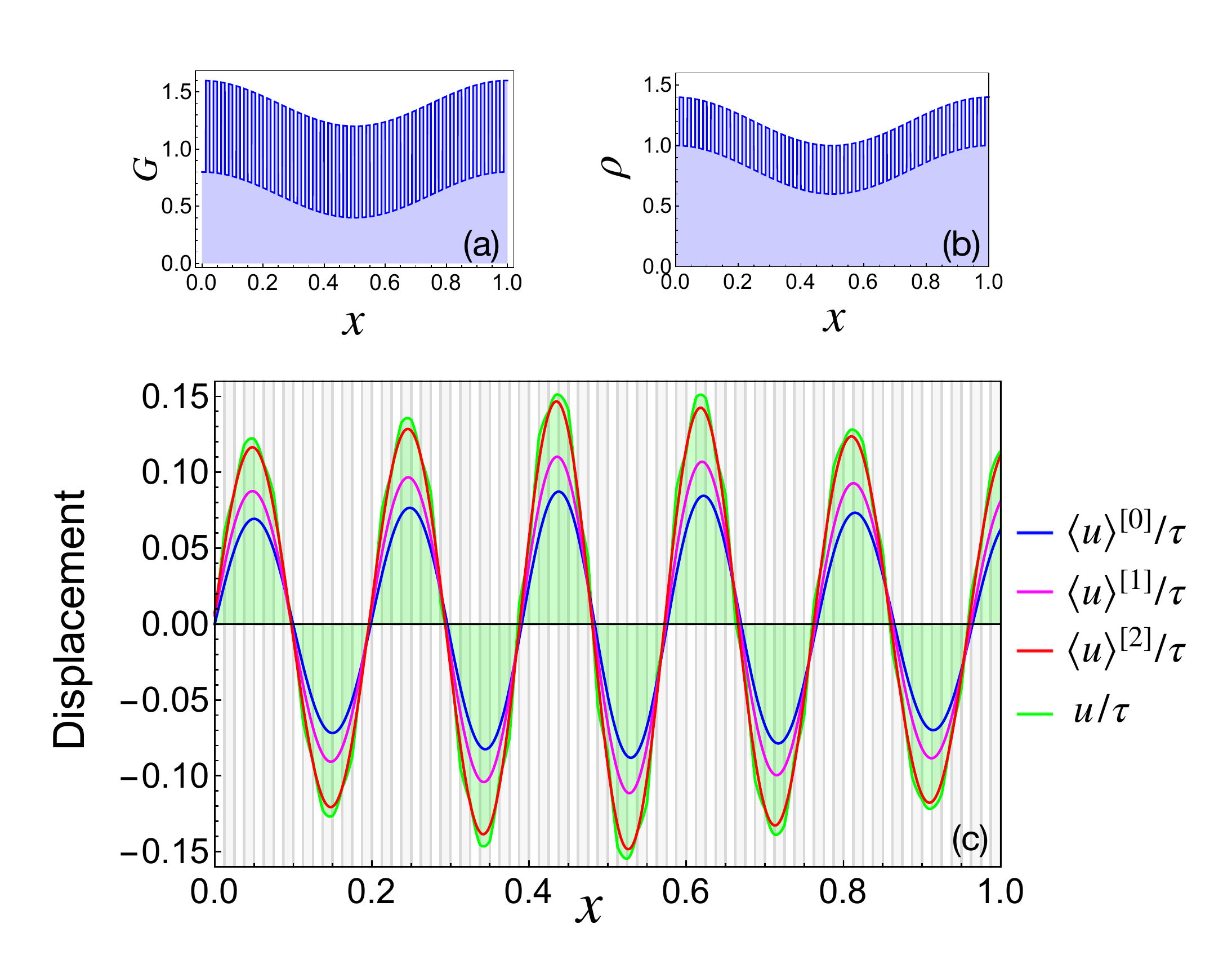}}\vspace*{-7mm}
\caption{Example Ex2: (a) Shear modulus profile, (b) mass density profile, and (c) ``exact'' wave motion $u$ versus homogenized mean fields $\langle u \rangle^{\mbox{\tiny{[0]}}}$, $\langle u \rangle^{\mbox{\tiny{[1]}}}$, and $\langle u \rangle^{\mbox{\tiny{[2]}}}$ for the BVP~\eqref{PB}.}
\label{PB_EX2_U}
\end{figure}

\begin{figure}[h!]  \vspace*{-1mm}
\centering{\includegraphics[width=0.77\linewidth]{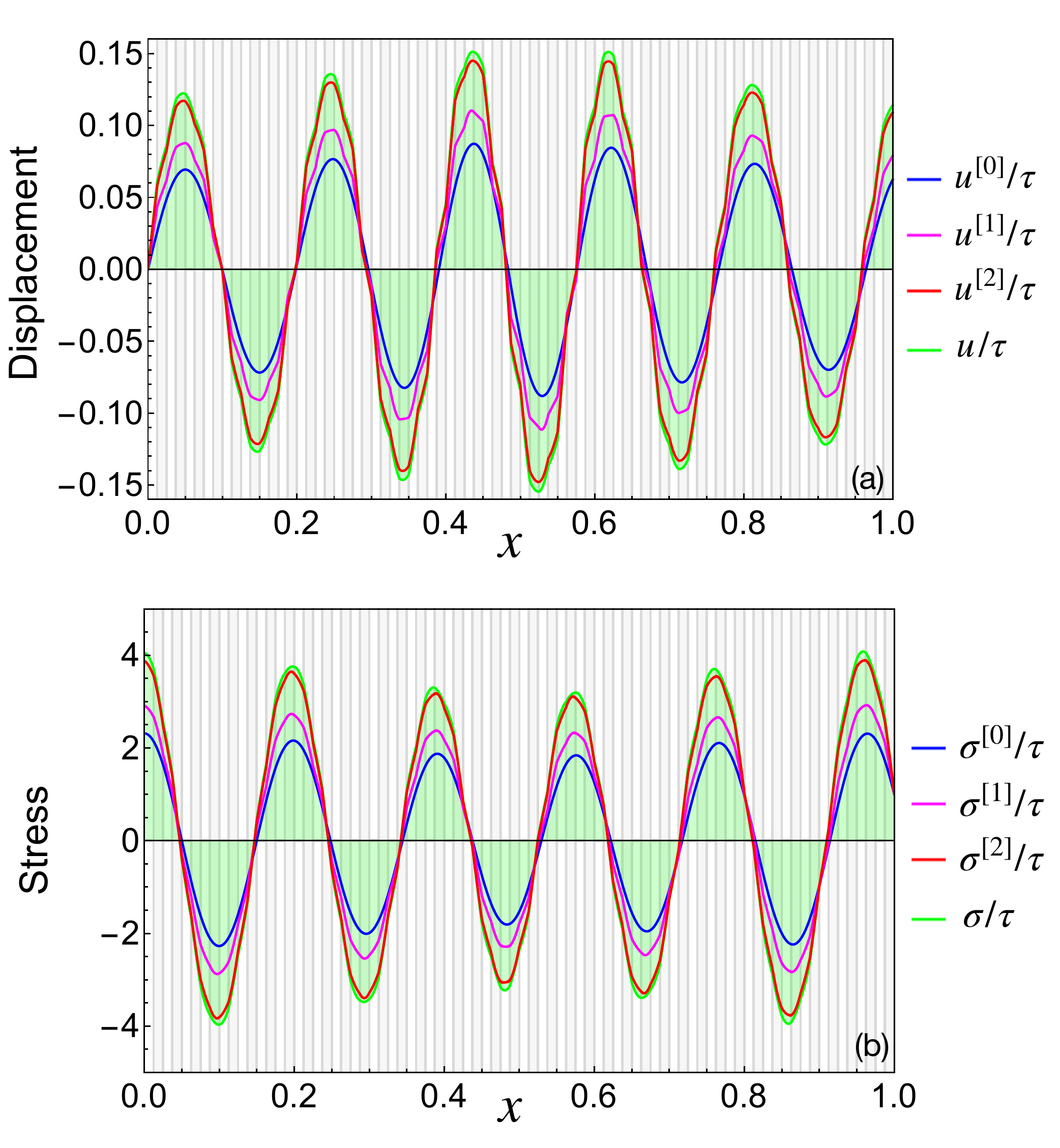}}\vspace*{-3mm}
\caption{Example Ex2: (a) ``exact'' wave motion $u$ versus homogenized approximations $u^{\mbox{\tiny{[0]}}}$, $u^{\mbox{\tiny{[1]}}}$, and $u^{\mbox{\tiny{[2]}}}$; (b) ``exact'' stress filed $\sigma(x)$ versus homogenized approximations $\sigma^{\mbox{\tiny{[0]}}}$, $\sigma^{\mbox{\tiny{[1]}}}$, and $\sigma^{\mbox{\tiny{[2]}}}$, for the BVP~\eqref{PB}.}	\label{PB_EX2_us}
\end{figure}

\begin{figure}[h!]  \vspace*{-1mm}
\centering{\includegraphics[width=0.8\linewidth]{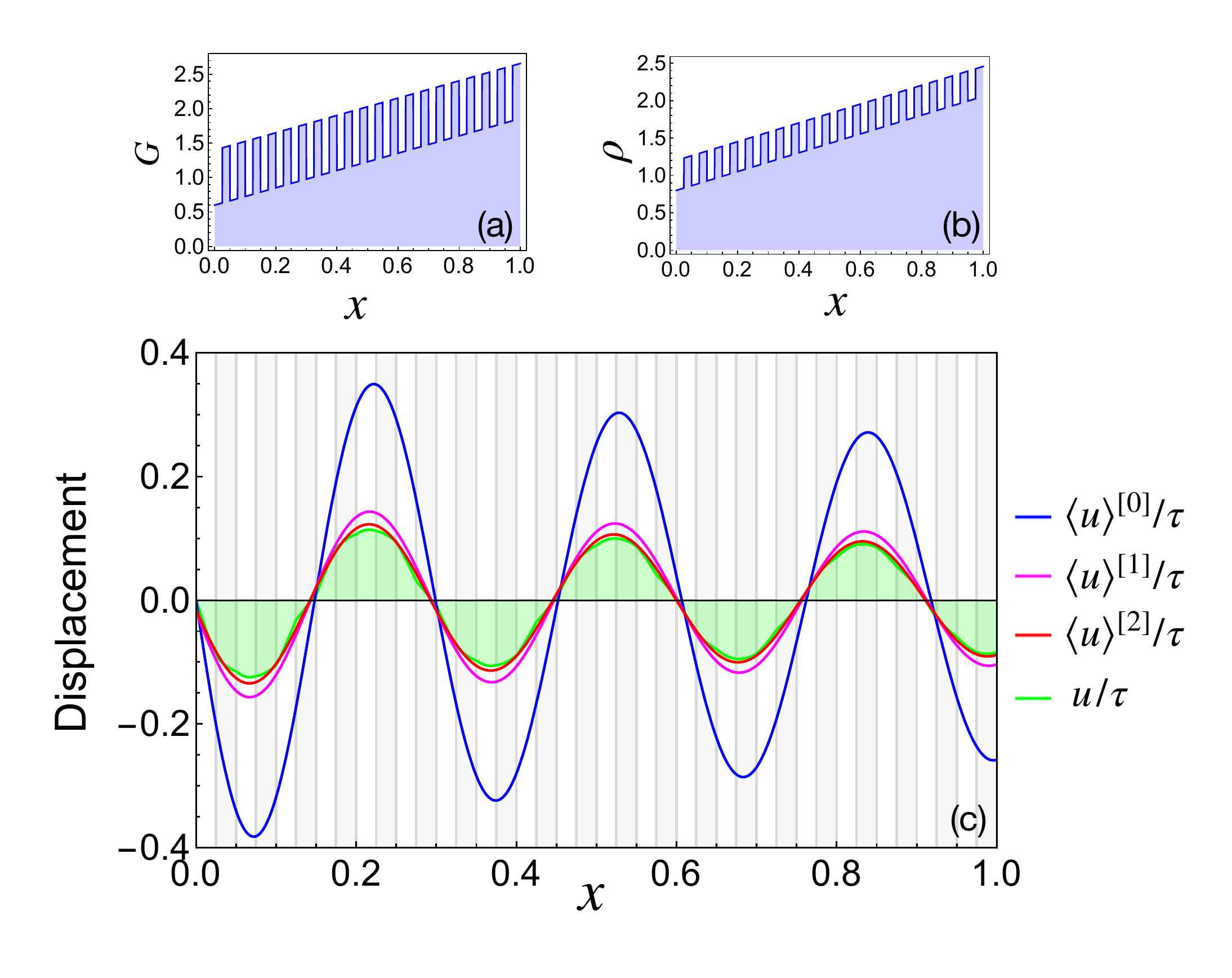}}\vspace*{-8mm}
\caption{Example Ex3: (a) Shear modulus profile, (b) mass density profile, and (c) ``exact'' wave motion $u$ versus homogenized mean fields $\langle u \rangle^{\mbox{\tiny{[0]}}}$, $\langle u \rangle^{\mbox{\tiny{[1]}}}$, and $\langle u \rangle^{\mbox{\tiny{[2]}}}$ for the BVP~\eqref{PB}}
\label{PB_EX3_U}
\end{figure}

\begin{figure}[h!]  \vspace*{-1mm}
\centering{\includegraphics[width=0.8\linewidth]{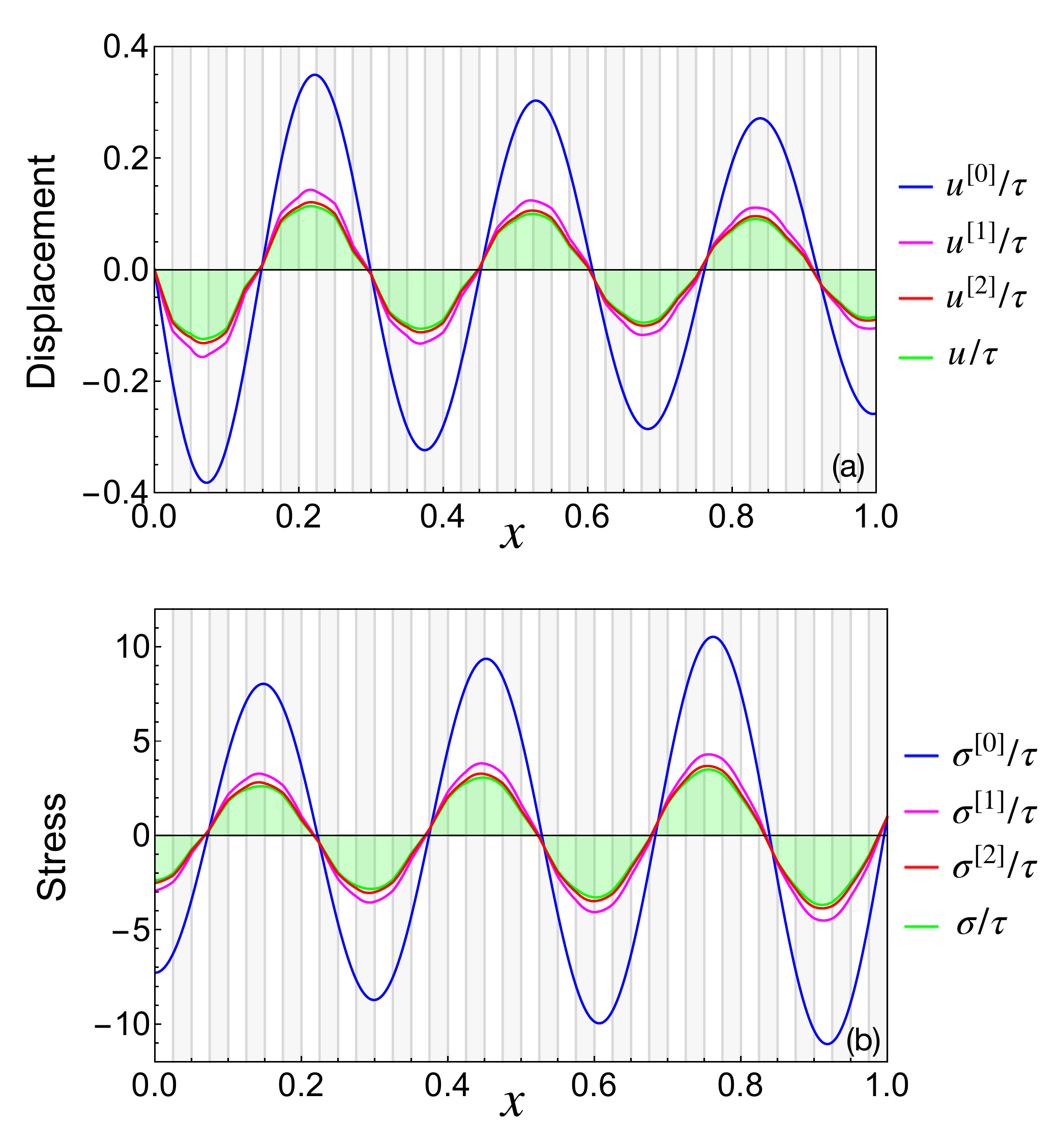}}\vspace*{-5mm}
\caption{Example Ex3: (a) ``exact'' wave motion $u(x)$ versus homogenized approximations $u^{\mbox{\tiny{[0]}}}$, $u^{\mbox{\tiny{[1]}}}$, and $u^{\mbox{\tiny{[2]}}}$; (b) ``exact'' stress filed $\sigma(x)$ versus homogenized approximations $\sigma^{\mbox{\tiny{[0]}}}$, $\sigma^{\mbox{\tiny{[1]}}}$, and $\sigma^{\mbox{\tiny{[2]}}}$, for the BVP~\eqref{PB}.}
\label{PB_EX3_us}
\end{figure}

To illustrate the performance of the homogenized models,  Fig.~\ref{PB_EX1_U}, Fig.~\ref{PB_EX2_U} and Fig.~\ref{PB_EX3_U} compare the \emph{mean fields} $\langle u \rangle^{\mbox{\tiny{[0]}}}$, $\langle u \rangle^{\mbox{\tiny{[1]}}}$, and $\langle u \rangle^{\mbox{\tiny{[2]}}}$ with the ``exact'' solution for examples Ex$1$, Ex$2$, and Ex$3$, respectively. As suggested earlier, the use of ($O(\eps)$ and~$O(\eps^2)$) asymptotic corrections is in this case critical to ensure the fidelity of the effective model. With such mean fields at hand, Fig.~\ref{PB_EX1_us}, Fig.~\ref{PB_EX2_us} and Fig.~\ref{PB_EX3_us} compare the \textit{full asymptotic approximations} $u^{\mbox{\tiny{[0]}}}$, $u^{\mbox{\tiny{[1]}}}$, and $u^{\mbox{\tiny{[2]}}}$ (computed via~\eqref{PB_u1} and~\eqref{PB_u2}) with the ``exact'' solution in terms of both displacement and stress waveforms, respectively, for examples Ex$1$, Ex$2$, and Ex$3$. A common observation from these displays is that the second-order model provides a satisfactory description of the exact wavefield, whereas its lower-order companions appear to be deficient. This contrast is especially striking in terms of leading-order model~$u^{\mbox{\tiny{[0]}}}=\langle u\rangle^{\mbox{\tiny{[0]}}}$ which appears to either undershoot by roughly~50\%, or overshoot by over~100\%, the ``exact'' solution. As mentioned earlier, however, all three approximations~($u^{\mbox{\tiny{[j]}}},\; j=\overline{0,2}$) are numerically observed to approach the ``exact'' solution as the excitation frequency~$\omega$ is gradually decreased relative to the values listed in Table~\ref{Table2}. 

\begin{remark}
By comparing the respective zero crossings in Figs.~\ref{PB_EX1_U}--\ref{PB_EX3_us}, it is apparent that even the lower-order models are quite good in capturing the phase of the solution -- a result that is consistent with the findings of Section~\ref{dispersion}. However it is also clear that, at least for the frequencies selected, lower-order approximations are inadequate for synthesizing the actual waveforms in quasi-periodic media. 
\end{remark}

\section{Summary} \label{Summary}

\noindent In this study, we pursue an effective description of the low-frequency wave motion in a macroscopically heterogeneous medium endowed with periodic microstructure. To this end, we deploy the framework of multiple scales and we apply the analysis to the scalar wave equation in one and multiple spatial dimensions. Through asymptotic expansion, the effective governing equation -- free of microscopic fluctuations -- is pursued up to the second order and shown to expose an intimate interplay between the dispersive effects of (periodic) micro-scale heterogeneities and their (generally non-periodic) macroscopic counterpart. More specifically, the germane low-frequency behavior is synthesized via a fourth-order differential equation (with smoothly varying coefficients) governing the mean wave motion in the medium, where the effect of microstructure is upscaled by way of the so-called cell functions. In an effort to demonstrate the relevance of our analysis toward solving boundary value problems, we also develop effective boundary conditions, up to the second order of asymptotic approximation, applicable to one-dimensional (1D) mean wave motion in a quasi-periodic medium. To our knowledge, this problem has escaped the scrutiny of earlier studies. We illustrate the analysis numerically in 1D by considering (i) low-frequency wave dispersion, (ii) mean-field homogenized description of waves propagating in a finite domain, and (iii) full-field homogenized description thereof. Specifically, we find that the microstructure may have a major effect on the overall wave dispersion in a quasi-periodic medium. In contrast to (i), however, where the latter appears to be well captured even by the leading-order model, the results in~(ii) and~(iii) illustrate the critical role that  higher-order corrections may have in maintaining the fidelity of homogenized waveform description.

\pagebreak

\section*{Acknowledgements} \label{Ack}

\noindent This work was supported in part through the endowed Shimizu Professorship, and Sommerfeld Fellowship to DS (Department of Civil, Environmental, and Geo- Engineering, University of Minnesota). The support provided by the Minnesota Supercomputing institute is kindly acknowledged. DS would also like to thank Othman Oudghiri-Idrissi for fruitful discussions and remarks.

\appendix\section{Cell functions describing effective wave motion in~$\mathbb{R}^d $ $(d>1)$} \label{appen1}

\noindent With reference to~\eqref{MicroS}, consider the situation where the unit cell~$Y$ is composed of~$N$ subdomains~$Y_q$ ($q=\overline{1,N}$) such that~$G''(\by)$ and~$\rho''(\by)$ according to either~\eqref{HetMat1} or~\eqref{HetMat2} vary smoothly within each~$Y_q$. In this setting, one finds that the (zero-mean) cell functions~$\boldsymbol{P}\in\mathbb{R}^d,\boldsymbol{Q}\in\mathbb{R}^{d\times d},\boldsymbol{R}\in\mathbb{R}^{d\times d\times d}, \tilde{\boldsymbol{P}}\in\mathbb{R}^d, \tilde{\boldsymbol{Q}}\in\mathbb{R}^{d\times d}$ and~$\tilde{\boldsymbol{R}}\in\mathbb{R}^{d}$ specifying the effective tensor coefficients in~\eqref{MMEq_2d} according to~\eqref{rho0_mu0_2d}, \eqref{xiC_2d} and~\eqref{phiC_2d} solve the respective boundary value problems 
\begin{eqnarray} \label{CF_P_2d}
\begin{split}
&\nabla_{\!\by } \sip  [ G(\boldsymbol{I}_2+ \nabla_{\!\by } \boldsymbol{P}) ] \;=\; \bzero,  \quad  \by \in  Y_q\\
&\boldsymbol{P}, ~ G\es \boldsymbol{n} \sip  (\boldsymbol{I}_2+ \nabla_{\!\by } \boldsymbol{P}), \quad Y\text{-periodic}; \quad \langle  \boldsymbol{P} \rangle = \bzero  \\
&  \llbracket \boldsymbol{P}\rrbracket =\bzero, \quad \llbracket G\es\boldsymbol{n} \sip (\boldsymbol{I}_2+ \nabla_{\!\by } \boldsymbol{P})\rrbracket =\bzero, \quad  \by \in \partial Y_q\backslash \partial Y;
\end{split}
\end{eqnarray}
\begin{eqnarray} \label{CF_Q_2d}
\begin{split}
& \nabla_{\!\by }  [G ( \boldsymbol{P} +\nabla_{\!\by } \! \cdot\! \boldsymbol{Q} ) ] \;=\; \frac{\rho}{\varrho^{\mbox{\tiny{(0)}}}} \, \boldsymbol{\mu^{\mbox{\tiny{(0)}}}} -G (\boldsymbol{I}_2 + \nabla_{\!\by } \boldsymbol{P}),  \quad \by \in  Y_q\\
&  \boldsymbol{Q}, ~ G\es \boldsymbol{n} \sip  ( \boldsymbol{P} +\nabla_{\!\by } \! \cdot\! \boldsymbol{Q} ), \quad Y\text{-periodic}; \quad \langle  \boldsymbol{Q} \rangle = \bzero  \\
&  \llbracket \boldsymbol{Q}\rrbracket =\bzero, \quad \llbracket  G\es\boldsymbol{n} \sip  ( \boldsymbol{P} +\nabla_{\!\by } \! \cdot\! \boldsymbol{Q} ) \rrbracket =\bzero,\quad   \by \in \partial Y_q\backslash \partial Y;
\end{split}
\end{eqnarray}
\begin{eqnarray} \label{CF_Ptilde_2d}
\begin{split}
& \nabla_{\!\by } \sip [G ( \nabla_{\!\bx} \boldsymbol{P} +\nabla_{\!\by } \tilde{\boldsymbol{P}} ) ] \;=\; \frac{\rho}{\varrho^{\mbox{\tiny{(0)}}} } \, \nabla \sip \boldsymbol{\mu^{\mbox{\tiny{(0)}}}}  - \nabla_{\!\bx} \sip [ G (\boldsymbol{I}_2 + \nabla_{\!\by } \boldsymbol{P}) ], \quad \by \in  Y_q\\
&  \tilde{\boldsymbol{P}}, ~ G\es\boldsymbol{n} \sip ( \nabla_{\!\bx} \boldsymbol{P} +\nabla_{\!\by } \tilde{\boldsymbol{P}} ), \quad Y\text{-periodic}; \quad \langle  \tilde{\boldsymbol{P}} \rangle = \bzero  \\
&  \llbracket \tilde{\boldsymbol{P}}\rrbracket =\bzero, \quad \llbracket  G\es\boldsymbol{n} \sip  ( \nabla_{\!\bx} \boldsymbol{P} +\nabla_{\!\by } \tilde{\boldsymbol{P}} ) \rrbracket =\bzero, \quad   \by \in \partial Y_q\backslash \partial Y;
\end{split}
\end{eqnarray}
\begin{eqnarray} \label{CF_Qtilde_2d}
\begin{split}
& \nabla_{\!\by } [G(\tilde{\boldsymbol{P}} + \nabla_{\!\bx}\sip\boldsymbol{Q} + \nabla_{\!\by}\sip\tilde{\boldsymbol{Q}})] \;=\; 
\frac{\rho}{\varrho^{\mbox{\tiny{(0)}}} } \big(\boldsymbol{\eta}+\nabla\sip\boldsymbol{\mu^{\mbox{\tiny{(1)}}}}\big)
- \nabla_{\!\bx} [G(\boldsymbol{P}+\nabla_{\!\by}\sip\boldsymbol{Q})] - G(\nabla_{\!\bx} \boldsymbol{P}+\nabla_{\!\by } \tilde{\boldsymbol{P}}) \\
& \hspace*{70mm} +  \rho \Big(\boldsymbol{P}-\frac{\boldsymbol{\varrho^{\mbox{\tiny{(1)}}}}}{\varrho^{\mbox{\tiny{(0)}}}}\Big) \otimes  \Big( \frac{\nabla \sip  \boldsymbol{\mu^{\mbox{\tiny{(0)}}}}}{\varrho^{\mbox{\tiny{(0)}}} }+\nabla \sip \big[\frac{\boldsymbol{\mu^{0}} }{\varrho^{\mbox{\tiny{(0)}}}}\big] \Big), \quad \by \in  Y_q \\
&  \tilde{\boldsymbol{Q}}, ~ G\es\boldsymbol{n} \sip (\tilde{\boldsymbol{P}} + \nabla_{\!\bx}\sip\boldsymbol{Q} + \nabla_{\!\by}\sip\tilde{\boldsymbol{Q}}), \quad Y\text{-periodic}; \quad \langle  \tilde{\boldsymbol{Q}} \rangle = \bzero  \\
&  \llbracket \tilde{\boldsymbol{Q}}\rrbracket =\bzero, \quad \llbracket G\es\boldsymbol{n} \sip (\tilde{\boldsymbol{P}} + \nabla_{\!\bx}\sip\boldsymbol{Q} + \nabla_{\!\by}\sip\tilde{\boldsymbol{Q}})\rrbracket =\bzero, \quad   \by \in \partial Y_q\backslash \partial Y;
\end{split}
\end{eqnarray}
\begin{eqnarray} \label{CF_R_2d}
\begin{split}
&  \nabla_{\!\by } [G(\boldsymbol{Q} + \nabla_{\!\by}\sip\boldsymbol{R})] \;=\; \frac{\rho}{\varrho^{\mbox{\tiny{(0)}}}} \boldsymbol{\mu^{\mbox{\tiny{(1)}}}}
- G(\boldsymbol{I}_2 \otimes \boldsymbol{P} + \nabla_{\!\by } \boldsymbol{Q})  
+ \rho \frac{\boldsymbol{\mu^{\mbox{\tiny{(0)}}}}}{\varrho^{\mbox{\tiny{(0)}}}} \otimes  \Big(\boldsymbol{P}-\frac{\boldsymbol{\varrho^{\mbox{\tiny{(1)}}}}}{\varrho^{\mbox{\tiny{(0)}}}}\Big), \quad \by \in  Y_q\\
&  \boldsymbol{R}, ~ G\es\boldsymbol{n} \sip (\boldsymbol{Q} + \nabla_{\!\by}\sip\boldsymbol{R}), \quad Y\text{-periodic}; \quad \langle  \boldsymbol{R} \rangle = \bzero  \\
&  \llbracket \boldsymbol{R}\rrbracket =\bzero, \quad \llbracket G\es\boldsymbol{n} \sip (\boldsymbol{Q} + \nabla_{\!\by}\sip\boldsymbol{R})\rrbracket  =\bzero, \quad   \by \in \partial Y_q\backslash \partial Y, 
\end{split} 
\end{eqnarray}
and 
\begin{eqnarray} \label{CF_Rtilde_2d}
\begin{split}
& \nabla_{\!\by } \sip  \big [G(\nabla_{\!\bx} \tilde{\boldsymbol{P}} + \nabla_{\!\by } \tilde{\boldsymbol{R}})]  \;=\; \frac{\rho}{\varrho^{\mbox{\tiny{(0)}}} } \nabla \sip  \boldsymbol{\eta } - \nabla_{\!\bx} \sip [G ( \nabla_{\!\bx} \boldsymbol{P} +\nabla_{\!\by } \tilde{\boldsymbol{P}} ) ] 
+ \rho \Big(\boldsymbol{P}-\frac{\boldsymbol{\varrho^{\mbox{\tiny{(1)}}}}}{\varrho^{\mbox{\tiny{(0)}}}}\Big) \es \nabla \sip \big[\frac{\nabla \sip  \boldsymbol{\mu^{0}}}{\varrho^{\mbox{\tiny{(0)}}}}\big], \quad \by \in  Y_q \\
& \tilde{\boldsymbol{R}}, ~ G\es\boldsymbol{n} \sip (\nabla_{\!\bx} \tilde{\boldsymbol{P}} + \nabla_{\!\by } \tilde{\boldsymbol{R}}), \quad Y\text{-periodic}; \quad \langle  \tilde{\boldsymbol{R}} \rangle = \bzero  \\
&  \llbracket \tilde{\boldsymbol{R}}\rrbracket =\bzero, \quad \llbracket G\es\boldsymbol{n} \sip (\nabla_{\!\bx} \tilde{\boldsymbol{P}} + \nabla_{\!\by } \tilde{\boldsymbol{R}})\rrbracket =\bzero, \quad   \by \in \partial Y_q\backslash \partial Y.
\end{split}
\end{eqnarray}

\bibliography{mybibfile} 

\end{document}